\documentclass{article}
\title{
Model selection for the robust efficient signal processing  observed with small L\'evy noise.
\thanks{
This work was done under financial support of the RSF Grant Number 14-49-00079 (National Research University “MPEI” 14 Krasnokazarmennaya, 111250 Moscow, Russia)
and by
 the RSF grant 17-11-01049 (National Research Tomsk State University).
}
}

\author{Slim Beltaief 
\thanks{
Laboratoire de Math\'ematiques Raphael Salem,
 UMR 6085 CNRS- Universit\'e de Rouen Normandie,  France,
e-mail: beltaiefslim@hotmail.fr
}
,
Oleg
Chernoyarov 
\thanks{
Department of Radio Engineering Devices and Antenna Systems,
National Research University “Moscow Power Engineering Institute”
and
ILSSP\&QF, National Research Tomsk State University,
 e-mail: chernoyarovov@mpei.ru
} 
 and
 Serguei
Pergamenchtchikov \thanks{
 Laboratoire de Math\'ematiques Raphael Salem,
 UMR 6085 CNRS- Universit\'e de Rouen Normandie,  France
and
ILSSP\&QF, National Research Tomsk State University,
 e-mail:
Serge.Pergamenchtchikov@univ-rouen.fr } }

\date{}

\usepackage{setspace}

\usepackage{epsfig}
\usepackage{graphicx}
\usepackage{wrapfig}
\usepackage{amssymb}
\usepackage{amsfonts}
\usepackage{amsmath}
\usepackage{amsthm}
\usepackage{enumerate}
\usepackage{multicol}

\newtheorem{theorem}{Theorem}[section]
\newtheorem{proposition}[theorem]{Proposition}
\newtheorem{lemma}[theorem]{Lemma}

\newtheorem{remark}{Remark}[section]
\newtheorem{corollary}[theorem]{Corollary}
\newcommand\fdem{$\Box$}
\newcommand\cA{{\cal A}}

\newcommand\cF{{\cal F}}

\newcommand\cB{{\cal B}}

\newcommand\cN{{\cal N}}

\newcommand\cX{{\cal X}}

\newcommand\cQ{{\cal Q}}
\newcommand\cS{{\cal S}}
\newcommand\cR{{\cal R}}

\newcommand\ve{\varepsilon}

\def\bbr{{\mathbb R}}

\def\text#1{\hbox{#1}}
\def\proof{{\noindent \bf Proof. }}
\def\endproof{\mbox{\ $\qed$}}
\def\card{\mbox{card}}

\def\E{{\bf E}}

\def\P{{\bf P}}
\def\C{{\bf C}}
\def\D{{\bf D}}

\def\U{{\bf U}}

\def\u{{\bf u}}
\def\p{{\bf p}}
\def\r{{\bf r}}
\def\c{{\bf c}}
\def\b{{\bf b}}
\def\l{{\bf l}}
\def\m{{\bf m}}
\def\L{{\bf L}}

\newcommand{\wh}{\widehat}
\newcommand{\wt}{\widetilde}

\newcommand\Er{\mbox{Err}}

\newcommand\Tr{\mbox{Tr}}
\def\R{{\bf R}}
\def\Chi{{\bf 1}}
\def\d{\mathrm{d}}
\def\build #1_#2{\mathrel{\mathop{\kern 0pt #1}\limits_\zs{#2}}}
\newcommand{\zs}[1]{{\mathchoice{#1}{#1}{\lower.25ex\hbox{$\scriptstyle#1$}}
{\lower0.25ex\hbox{$\scriptscriptstyle#1$}}}}
\numberwithin{equation}{section}

\begin{document}

\maketitle

\begin{abstract}

We develop a new model selection method for  the adaptive robust efficient nonparametric signal estimation observed   
with impulse noise which is defined by the general non Gaussian L\'evy processes.  On the basis of the developed method, we construct the estimation procedures which are analyzed in two settings: in non asymptotic and asymptotic ones.  For the first time for such models  we show non asymptotic sharp oracle inequalities for the quadratic and for the robust risks,  i.e.  we show that  the constructed procedures are optimal in the  sharp oracle inequalities sense. 
Next, by making use of the obtained oracle inequalities, we provide the asymptotic efficiency property for the developed estimation methods in the adaptive setting when the  signal/noise ratio goes to  infinity.  We apply the developed model selection methods for the signals number detection problem in  multi-path  information transmission.
 
 \end{abstract}

\vspace*{5mm}

\vspace*{5mm}

\newpage

\section{Introduction}\label{sec:In}

In this paper we consider the
signal estimation problem
on the basis of  observations defined by
  the nonparametric regression model in continuous time
with pulse noises of small intensity, i.e., 
\begin{equation}\label{sec:In.1}
 \d\,y_\zs{t} = S(t)\d\,t + \varepsilon\,
\d \xi_\zs{t}\,,\quad 
 0\le t \le 1\,,
\end{equation}
where $S(\cdot)$ is an unknown deterministic signal (i.e., $[0,1]\to\bbr$ nonrandom function), 
$(\xi_\zs{t})_\zs{0\le t\le 1}$ is an unobserved noise and
$\varepsilon>0$ is the noise intensity. The problem is to estimate
the function $S$ on the observations $(y_\zs{t})_\zs{0\le t\le 1}$ when $\varepsilon\to 0$.
 Note that if $(\xi_\zs{t})_\zs{0\le t\le 1}$ is a brownian motion, then we obtain the  "signal+white noise" model
  which is very popular in statistical radio-physics and
 is well studied by many authors: Ibragimov and Khasminskii in \cite{IbragimovKhasminskii1981}, Pinsker in \cite{Pinsker1981}, 
Kutoyants in \cite{Kutoyants1984} and \cite{Kutoyants1994},  etc..
 The condition $\varepsilon\to 0$ means that the signal/noise ratio goes to  infinity.
 In this paper,  we assume that in addition 
 to the intrinsic noise in the radio-electronic system, approximated usually by the gaussian white noise, 
 the useful signal $S$ is distorted by the impulse noise flow defined
by  L\'evy process with  jumps introduced in the next section.
The cause of the appearance of the pulse stream in the radio-electronic systems can be, for example, either  external unintended (atmospheric) noises,  intentional impulse noises or  errors in the demodulation and  channel decoding for the binary information symbols. 
Note that,  the impulse noises
for the signal detection  problems
 have been introduced  for the first time  by Kassam in
\cite{Kassam1988} on the basis of the compound Poisson processes. 
 Later, 
 Konev, Pergamenshchikov and Pchelintsev used
 the
  compound Poisson processes
  in
  \cite{Pchelintsev2013, KonevPergamenshchikovPchelintsev2014}
for the parametric regression models
and
   in 
  \cite{KonevPergamenshchikov2012, KonevPergamenshchikov2015}
for the nonparametric signal estimation problems.
 However,  the compound Poisson process can  describe only the large impulses  influence of small frequencies.
 It should be noted that in the telecommunication systems,  the noise impulses are without limitations on frequencies and therefore,
 the compound Poisson models are too restricted for  practical applications.   
 To include  all possible impulse noises, we propose to use a general non-gaussian   L\'evy  processes  
 in 
 the observation model  \eqref{sec:In.1}.
 In this paper,  we consider a nonparametric estimation problem in the adaptive setting,
i.e., when the regularity of the signal $S$ is unknown. Moreover,  we also assume that the
distribution $Q$
of
 the noise process  $(\xi_\zs{t})_\zs{0\le t\le 1}$  is  unknown.  
 It is only known  that this distribution belongs to the
  distribution family $\cQ^{*}_\zs{\varepsilon}$ defined  in the next section.
    By these reasons,  we use  the robust estimation approach proposed for nonparametric problems by
    Galtchouk, Konev and Pergamenshchikov
    in
 \cite{GaltchoukPergamenshchikov2006, KonevPergamenshchikov2012, KonevPergamenshchikov2015}.
We set the robust risks as
\begin{equation}\label{sec:robust-risks}
\cR^{*}_\zs{\varepsilon}(\wh{S}_\zs{\varepsilon},S)=\sup_\zs{Q\in\cQ^{*}_\zs{\varepsilon}}\,
\cR_\zs{Q}(\wh{S}_\zs{\varepsilon},S) \,,
\end{equation}
where  $\wh{S}_\zs{\varepsilon}$ is an estimator (i.e., any measurable function of $(y_\zs{t})_\zs{0\le t\le 1}$),
\begin{equation}\label{sec:risks}
\cR_\zs{Q}(\wh{S}_\zs{\varepsilon},S):=
\E_\zs{Q,S}\,\|\wh{S}_\zs{\varepsilon}-S\|^2
\quad\mbox{and}\quad
\Vert S\Vert^{2}=\int^{1}_\zs{0}\,S^{2}(t)\d t
\,.
\end{equation}
In this paper, we develop a sharp model selection method for  estimating the unknown signal $S$.  
The interest to such statistical procedures can be explained by the fact that 
they provide
adaptive solutions for the nonparametric estimation through the non-asymptotic
oracle inequalities which give the  non-asymptotic upper
bound for the quadratic risk including
the minimal risk over chosen  family of estimators
with some coefficient close to one. Such inequalities were obtained, for example, 
by 
Galtchouk and Pergamenshchikov 
\cite{GaltchoukPergamenshchikov2009a} for 
 non Gaussian regression models in
discrete time
and by Konev and Pergamenshchikov 
\cite{KonevPergamenshchikov2009a} for general regression semimartingale models in continuous time. 
It should be noted that for the first time the model selection methods 
were proposed by 
 Akaike \cite{Akaike1974} and Mallows \cite{Mallows1973} 
for parametric models. Then, 
by using the  oracle inequalities approach,
these methods had been developed
for the nonparametric estimation by
Barron, Birg\'e, Massart \cite{BarronBirgeMassart1999},  
for  Gaussian regression models
and
by Fourdrinier and Pergamenshchikov \cite{FourdrinierPergamenshchikov20007}
for  non Gaussian models. 
We know that an oracle inequality yields the upper bound for the  risks via  minimal
risk corresponding to a chosen estimators family.
Unfortunately, the oracle inequalities obtained in these papers
 can not be used
for
 the efficient estimation in the adaptive setting, since
the upper bounds in these inequalities
have some fixed coefficients in the main terms which are more than one. 
In order to provide the efficiency property for  model selection
procedures, one needs to obtain the sharp oracle inequalities, i.e., 
 in which the coefficient at the principal term on the right-hand side of the inequality is close to one.
To obtain  such inequalities for  general  non gaussian observations, one needs to use the model selection method
based on the weighted least square estimators
 proposed by
 Galtchouk and Pergamenshchikov
 \cite{GaltchoukPergamenshchikov2009a, GaltchoukPergamenshchikov2009b}  for the heteroscedastic
 regression models in discrete time and developed then by 
Konev and Pergamenshchikov in
\cite{KonevPergamenshchikov2009a, KonevPergamenshchikov2009b, KonevPergamenshchikov2012, KonevPergamenshchikov2015}
for  semimartingale models in  continuous time,  i.e.,  when the observation process is given by the following stochastic differential
equation
\begin{equation}\label{sec:In.1.check_00}
  \d x_\zs{t}=S(t)\d t+
  \d \eta_\zs{t}\,,
  \quad 0\le t\le n\,,\quad (n\to\infty)\,,
 \end{equation}
where
$S$
is an unknown $1$ - 
periodic signal
and the unobserved noise
 $(\eta_\zs{t})_\zs{t\ge 0}$ is square integrated  semi-martingale.
Note that, for any $0<t<1$,  
setting
$\check{x}_\zs{t}=n^{-1}\,\sum^{n}_\zs{j=1}(x_\zs{t+j}-x_\zs{j})$,
we can represent this model as a model with  small parameter of form \eqref{sec:In.1}
\begin{equation}\label{sec:In.1.check}
\d \check{x}_\zs{t}=S(t)\d t+\varepsilon\,\d \check{\eta}_\zs{t}\,,
\end{equation}
where
$\varepsilon=n^{-1/2}$ and
$\check{\eta}_\zs{t}=n^{-1/2}\,\sum^{n-1}_\zs{j=0}(\eta_\zs{t+j}-\eta_\zs{j})$. 
If $(\eta_\zs{t})_\zs{t\ge 0}$ is   L\'evy process, then 
$\check{\eta}_\zs{t}$ is  L\'evy process as well. But 
the  main difference between the
models  
\eqref{sec:In.1}
and
\eqref{sec:In.1.check}
 is that
the jumps in the last one are small, i.e.,
\begin{equation}\label{sec:In.2.check}
\Delta \check{\eta}_\zs{t}
=
\check{\eta}_\zs{t}
-
\check{\eta}_\zs{t-}
 = \mbox{O}(n^{-1/2})= \mbox{O}(\varepsilon)
\quad\mbox{as}\quad
\varepsilon\to 0
\,.
\end{equation}
But there is  no such property in  model \eqref{sec:In.1}. It should be noted that  property \eqref{sec:In.2.check}
is crucial in the non asymptotic analysis  for  observations on  large time intervals, i.e.  the methods developed  for  model 
\eqref{sec:In.1.check_00}  can not be used for the  problem \eqref{sec:In.1}. Moreover, it should be emphasized that
the selection model methods proposed by Konev and Pergameshchikov for the model  \eqref{sec:In.1.check_00}
provide the adaptive efficient estimation  only for the case when the L\'evy measure is finite. 
This condition considerably reduces their  applications in practical problems.
So,  the main goal of this paper is to develop 
a new model selection method for  the adaptive efficient signal estimation problem
in a nonparametric regression \eqref{sec:In.1} 
for the general L\'evy noises without limitations on the jumps.
First, we construct some model selection procedures and we show the sharp non asymptotic oracle inequalities
for the  risks \eqref{sec:robust-risks} and \eqref{sec:risks}. 
To do this  in Proposition \ref{Pr.sec: MaPr.1} we develop a special analytical tool
to study the non asymptotic behavior of the jumps in the model  \eqref{sec:In.1}
with the infinite (or finite) L\'evy measure. 
Moreover, to study the efficiency we develop the  Van Trees method 
for general L\'evy processes 
and we obtain in Proposition \ref{Pr.sec:VanTreesIn++} a new lower bound for quadratic risks in the model \eqref{sec:In.1}.
Then, by making use of this lower bound
  we find the Pinsker constant. As to the upper bound,
 similarly to Konev and Pergamenshchikov
\cite{KonevPergamenshchikov2009b},
  we use the obtained  sharp oracle inequality for  
 weight least square estimators  containing the efficient Pinsker procedure. Therefore, through  oracle inequality,
 we estimate from above the risk for the constructed model selection procedure by the efficient risk 
 up to some coefficient which goes to one. As a result, we provide 
the robust efficiency property for the constructed procedure in  adaptive setting.
  As an application for the developed model selection method, in this paper,
 we consider the signals number detection problem for  model \eqref{sec:In.1}.
  In many areas of science and technology, this problem arises how to select the number of freedom degrees
  for a statistical model that most adequately describes phenomena under studies (see, for example, Akaike \cite{Akaike1974}).
  An important class of such problems is the detection problem of  signals  number  with unknown parameters observed  in  multi-path
  information transmission with  noises.  For example,  in the signal multi-path information transmission, there is a detection problem for the number of rays 
  in the multi-path channel.  This problem is often reduced to the detection of the number of signals. 
  As a result, the effective detection algorithms  can significantly improve the noise immunity in  data transmission over a multi-path channel
 (see, for example, the papers of  Flaksman,  Manelis, El-BeMac, Trifonov, Kharin, and Chernoyarov 
 \cite{Flaksman2002, Manelis2007, El-BeMac1978, TrifonovKharin2013, TrifonovShinakov1986, TrifonovKharinChernoyarovKalashnikov2015, TrifonovKharin2015_a, TrifonovKharin2015_b}). 
 These problems for signals with unknown amplitudes are discussed by Trifonov and Kharin in \cite{TrifonovKharin2013}. 
The signal amplitude is an energy parameter because it affects the signal energy. 
At the same time, quite often, such as in radars studied by El-BeMac in \cite{El-BeMac1978}, 
it is necessary to detect the number of   signals, which besides an unknown amplitude, 
 it contains non-energy parameters such as the time and directions of the  signal
arrival, its frequency and  initial phase. Moreover,  Trifonov, Kharin, Chernoyarov and Kalashnikov
 in \cite{TrifonovKharinChernoyarovKalashnikov2015}  considered this problem   with unknown initial phases, 
 in \cite{TrifonovKharin2015_a} with unknown amplitudes and phases and in \cite{TrifonovKharin2015_b} the 
 detection signal number problem is considered for orthogonal signals with arbitrary non-energy parameters. 
  In all these papers the signals number detection problems are considered only for observation with  Gaussian white nose.
In this paper we consider this problem for the non Gaussian impulse noise.

The rest of the paper is organized as follows.
In Section \ref{sec:Mconds}, we give the main conditions which will be assumed
 for the model \eqref{sec:In.1}.  In Section \ref{sec:Tr}, we transform the observation model to delete large jumps
and 
we develop an  analytical tool  for the L\'evy  regression models in continuous time which provides to 
 study the non asymptotic behavior for the 
 sum of the deviations of the squares of the stochastic integrals of basic functions with respect to the non Gaussian L\'evy processes.
In Section \ref{sec:Mo}, we construct the sharp model selection procedure. In Section \ref{sec:Mrs} we give the main results on the sharp oracle inequalities
and on the asymptotic robust efficiency. 
In Section \ref{sec: VanTreesIn__} we obtain the van Trees inequlaity for the general L\'evy processes.
 In Sections \ref{sec:Lobn} and  \ref{sec:UpBn}, we study the lower
and upper bounds  for the robust risks. In Section \ref{sec:NumSgn}, we study the signals number detection problem 
through the developed model selection method. In Section \ref{sec:Siml} we give  simulations results. Section \ref{sec:Pr}  contains
the proofs of all main results. In Appendix, we bring all   auxiliary results.

\bigskip

\section{Main conditions}\label{sec:Mconds}

In this section
  we assume that
 the noise process $(\xi_\zs{t})_\zs{0\le t\le 1}$
is defined as
\begin{equation}\label{sec:In.1+1}
\xi_\zs{t}=\varrho_\zs{1} w_\zs{t} + \varrho_\zs{2} z_\zs{t}
\quad\mbox{and}\quad
z_\zs{t}=x*(\mu-\wt{\mu})_\zs{t}
\,,
\end{equation}
where, $\varrho_\zs{1}$ and $\varrho_\zs{2}$ are  some unknown constants,
$(w_\zs{t})_\zs{0\le t\le\,1}$ is a standard brownian motion, 
 "*" denotes the stochastic integral with respect to the compensated  jump measure (see, for example in
 Jacod and Shiryaev
\cite{JacodShiryaev2002} or Cont and Tankov \cite{ContTankov2004} for details), 
 $\mu(\d s\,\d x)$ is a jump measure with  deterministic
compensator $\wt{\mu}(\d s\,\d x)=\d s\Pi(\d x)$, 
$\Pi(\cdot)$ is the unknown L\'evy measure, i.e.  some positive measure on $\bbr_\zs{*}=\bbr\setminus \{0\}$, 
such that
\begin{equation}\label{sec:Ex.1-00}
\Pi(x^{2})=1
\quad\mbox{and}\quad
\Pi(x^{4})
\,<\,\infty\,,
\end{equation}
where $\Pi(\vert x\vert^{m})=\int_\zs{\bbr_\zs{*}}\,\vert z\vert^{m}\,\Pi(\d z)$. Note that
 the measure $\Pi(\bbr_\zs{*})$ could be equal to $+\infty$.
In the sequel we will denote by $Q$ the distribution of the process $(\xi_\zs{t})_\zs{0\le t\le 1}$.
We assume that the parameters
    $\varrho_\zs{1}$
and $\varrho_\zs{2}$ satisfy the conditions
\begin{equation}\label{sec:Ex.01-1}
0< \check{\varrho}_\zs{\varepsilon}\le \varrho^{2}_\zs{1}
\quad\mbox{and}\quad
\varkappa_\zs{Q}=\varrho^{2}_\zs{1}+\varrho^{2}_\zs{2}\,
\le 
\varsigma^{*}_\zs{\varepsilon}
\,,
\end{equation}
where
 the bounds
$\check{\varrho}_\zs{\varepsilon}$ and 
$\varsigma^{*}_\zs{\varepsilon}$ are such that for any $\b>0$
\begin{equation}\label{sec:Ex.01-2}
\liminf_\zs{\varepsilon\to 0}\,
\varepsilon^{-\b}\,
\,\check{\varrho}_\zs{\varepsilon}
>0
\quad\mbox{and}\quad
\lim_\zs{\varepsilon\to 0}\,\varepsilon^{\b}\,\varsigma^{*}_\zs{\varepsilon}
=0
\,.
\end{equation}
We denote
by $\cQ^{*}_\zs{\varepsilon}$ the family of all distributions $Q$ of the process \eqref{sec:In.1+1}
 in the Skorokhod space $\D[0,1]$ for which the conditions  \eqref{sec:Ex.01-1} and
\eqref{sec:Ex.01-2}
hold.

\bigskip

\section{Transformation of the observations}\label{sec:Tr}

First of all, we need to eliminate the large jumps in the observations
\eqref{sec:In.1}, i.e.  we transform this model as
\begin{equation}\label{sec:Mo.0}
\check{y}_\zs{t}=y_\zs{t}-\sum_\zs{0\le s\le t}\,\Delta y_\zs{s}\,\Chi_\zs{\{\vert \Delta y_\zs{s}\vert>\overline{a}\}}
\,.
\end{equation}
The parameter $\overline{a}=\overline{a}_\zs{\varepsilon}>0$ will be chosen later. So, we obtain that
\begin{equation}\label{sec:Mo.1}
\d \check{y}_\zs{t}=S(t)\d t+\varepsilon\d \check{\xi}_\zs{t}
- \varepsilon\,\varrho_\zs{2}\,\Pi(\overline{h}_\zs{\varepsilon})\,\d t
\,,
\end{equation}
where 
$
\check{\xi}_\zs{t}
=\varrho_\zs{1} w_\zs{t} + \varrho_\zs{2} \,\check{z}_\zs{t}$
and
$\check{z}_\zs{t}=h_\zs{\varepsilon}*(\mu-\wt{\mu})_\zs{t}$.
The functions
 $h_\zs{\varepsilon}(x)=x\Chi_\zs{\{\vert x\vert\le \check{\upsilon}_\zs{\varepsilon}\}}$ and $\overline{h}_\zs{\varepsilon}(x)=x\Chi_\zs{\{\vert x\vert> \check{\upsilon}_\zs{\varepsilon}\}}$ with the
 truncation threshold $\check{\upsilon}_\zs{\varepsilon}=\overline{a}/\vert\varrho_\zs{2}\vert\varepsilon$.

\begin{remark}\label{Re.sec:Ma.1++}
It should be noted that the sum in the transformation \eqref{sec:Mo.0}
is finite since the cadlag process has only finite number of jumps more than 
some positive threshold in absolute value.
\end{remark}

Let $(\phi_\zs{j})_\zs{j\ge\, 1}$ be an orthonormal  basis in $\L_\zs{2}[0,1]$ with $\phi_\zs{1}\equiv 1$.
We assume that this basis is uniformly bounded, i.e.
for some  constant $\phi^{*}> 0$,  which may be dependent on $\varepsilon>0$,
\begin{equation}\label{sec:In.3-00}
\sup_\zs{0\le j\le n}\,\sup_\zs{0\le t\le 1}\vert\phi_\zs{j}(t)\vert\,
\le\,
\phi^{*}
<\infty\,,
\end{equation}
where $n=n_\zs{\varepsilon}=[1/\varepsilon^{2}]$  and  $[x]$ denotes integer part of $x$. For example, we can take 
 the trigonometric basis    defined as $\Tr_\zs{1}\equiv 1$ and for $j\ge 2$
\begin{equation}\label{sec:In.5}
 \Tr_\zs{j}(x)= \sqrt 2
\left\{
\begin{array}{c}
\cos(2\pi[j/2] x)\, \quad\mbox{for even}\quad j \,;\\[4mm]
\sin(2\pi[j/2] x)\quad\mbox{for odd}\quad j\,.
\end{array}
\right.
\end{equation}

\noindent  Moreover, note that
 for any $[0,1] \to\bbr$ function $f$ from $\L_\zs{2}[0,1]$ and 
 for any $0\le t\le 1$
 the integrals
\begin{equation}\label{sec:In.2}
I_\zs{t}(f)=\int_\zs{0}^{t} f(s) \d\xi_\zs{s}
\quad\mbox{and}\quad
\check{I}_\zs{t}(f)=\int_\zs{0}^{t} f(s) \d\check{\xi}_\zs{s}
\end{equation}
are well defined with $\E\,I_\zs{t}(f)=0$, $\E\,\check{I}_\zs{t}(f)=0$, 
\begin{equation}\label{sec:In.3}
\E\,I^{2}_\zs{t}(f) =\varkappa_\zs{Q}\,\Vert f\Vert^{2}_\zs{t}
\quad\mbox{and}\quad
\E\,\check{I}^{2}_\zs{t}(f) =\check{\varkappa}_\zs{Q}\,\Vert f\Vert^{2}_\zs{t}
\,,
\end{equation}
where $\Vert f\Vert^{2}_\zs{t}=\int^{t}_\zs{0}\,f^{2}(s)\d s$
 and $\check{\varkappa}_\zs{Q}=\varrho^{2}_\zs{1}+\varrho^{2}_\zs{2}\Pi(h^{2}_\zs{\varepsilon})$. 
In the sequel we denote by 
$$
(f,g)_\zs{t}=\int^{t}_\zs{0}\,f(s)g(s)\,\d s
\quad\mbox{and}\quad 
(f,g)=\int^{1}_\zs{0}\,f(s)g(s)\,\d s
\,.
$$
To estimate the function $S$ we use the following Fourier series
\begin{equation}\label{sec:In.6}
S(t)=\sum_\zs{j\ge 1}\,\theta_\zs{j}\,\phi_\zs{j}(t)
\quad\mbox{and}\quad
\theta_\zs{j}=(S,\phi_\zs{j})\,.
\end{equation}
These coefficients can be estimated by the following way. The first we estimate as
$$
\wh{\theta}_\zs{1,\varepsilon}=  
\int_\zs{0}^{1}  \phi_\zs{1}(t) \d\,y_\zs{t}
=\theta_\zs{1}
+
\varepsilon \xi_\zs{1}
$$
and for $j\ge 2$
\begin{equation}\label{sec:In.7}
\wh{\theta}_\zs{j,\varepsilon}=  
\int_\zs{0}^{1}  \phi_\zs{j}(t) \d\,\check{y}_\zs{t}\,.
\end{equation}
Taking into account here that for $j\ge 2$ the integral $\int^{1}_\zs{0}\phi_\zs{j}(t)\d t=0$ we obtain
from \eqref{sec:Mo.1} that these Fourier coefficients can be represented as
$$
\wh{\theta}_\zs{j,\varepsilon}= \theta_\zs{j} + \varepsilon\, \overline{\xi}_\zs{j}
\quad\mbox{and}\quad
\overline{\xi}_\zs{j}=  \check{I}_\zs{1}(\phi_\zs{j})
\,.
$$
Setting $\overline{\xi}_\zs{1}=\xi_\zs{1}$ we obtain that for any $j\ge 1$
\begin{equation}\label{sec:In.8}
\wh{\theta}_\zs{j,\varepsilon}= \theta_\zs{j} + \varepsilon\, \overline{\xi}_\zs{j}
\,.
\end{equation}

Now, according to the model selection approach developed in 
\cite{KonevPergamenshchikov2009a, KonevPergamenshchikov2009b}
we need to study for any $u\in \bbr^{n}$ the following functions

\begin{equation}\label{sec: MaPr.1}
B_\zs{1,\varepsilon}(u)=\sum^{n}_\zs{j=1}\,u_\zs{j}\,
(\E_\zs{Q}\,\overline{\xi}^{2}_\zs{j}-\check{\varkappa}_\zs{Q})
\quad\mbox{and}\quad
B_\zs{2,\varepsilon}(u)=
\sum^{n}_\zs{j=1}\,u_\zs{j}\,\wt{\xi}_\zs{j}
\,,
\end{equation}
where $\wt{\xi}_\zs{j}= \overline{\xi}^{2}_\zs{j}- \E_\zs{Q}\,\overline{\xi}^{2}_\zs{j}$.

\begin{proposition}
\label{Pr.sec: MaPr.0}
The following upper bound holds
\begin{equation}
\label{sec: MaPr.1+1}
\sup_\zs{u\in[0,1]^{n}}\left\vert B_\zs{1,\varepsilon}(u) \right\vert
\le
\varkappa_\zs{Q}
\,.
\end{equation}
\end{proposition}
\proof
Note that $
\vert\E_\zs{Q}\,\overline{\xi}^{2}_\zs{1}-\check{\varkappa}_\zs{Q}\vert
=
\vert\E_\zs{Q}\,\xi^{2}_\zs{1}-\check{\varkappa}_\zs{Q}\vert
=
\varkappa_\zs{Q}
-\check{\varkappa}_\zs{Q}
\le  \varkappa_\zs{Q}$
 and 
$\E_\zs{Q}\,\overline{\xi}^{2}_\zs{j}=\check{\varkappa}_\zs{Q}$
for $j\ge 2$. So, from this we immediately obtain the upper bound \eqref{sec: MaPr.1+1}.
\endproof

\noindent 
  Now,  for any $u\in\bbr^{n}$ we set
\begin{equation}
\label{sec: MaPr.2--}
\vert u\vert^{2}=\sum^{n}_\zs{j=1}\,u^{2}_\zs{j}
\quad\mbox{and}\quad
\#(u)=\sum^{n}_\zs{j=1}\,\Chi_\zs{\{u_\zs{j}\neq 0\}}
\,.
\end{equation}

\begin{proposition}\label{Pr.sec: MaPr.1}
For any fixed truncation parameter $\overline{a}>0$ 
and for any vector  $u\in\bbr^{n}$ with $\vert u\vert\le 1$
\begin{equation}
\label{sec: MaPr.2}
\E_\zs{Q}\,
 B^{2}_\zs{2,\varepsilon}(u) 
\le U_\zs{Q}
+
6\check{\varkappa}_\zs{Q}\,
\left(
\frac{\overline{a}}{\varepsilon}
\right)^{2}
\,\#(u)\,(\phi^{*})^{4}
\,,
\end{equation}
where  $U_\zs{Q}=24 \varkappa^{2}_\zs{Q}+6\varrho^{4}_\zs{2}\,\Pi(x^{4})$.
\end{proposition}

\begin{remark}\label{Re.sec:Ma.2++1}
It should be noted that the last term in the non asymptotic  upper bound 
\eqref{sec: MaPr.2} is appeared due to the influence of the jumps in the observations
\eqref{sec:In.1}.
Note that we will use  the upper bounds 
\eqref{sec: MaPr.1+1} -- \eqref{sec: MaPr.2}
to obtain the non asymptotic sharp oracle oracle inequalities.
\end{remark}

\section{Model selection}\label{sec:Mo}

We estimate the function $S(x)$ for $x\in [0,1]$ by the weighted least squares estimator
\begin{equation}\label{sec:Mo.1++1}
\wh{S}_\lambda (x) = \sum_\zs{j=1}^{n} \lambda(j) \wh{\theta}_\zs{j,\varepsilon} \phi_\zs{j}(x)\,,
\end{equation}
where $n=[1/\varepsilon^2 ]$, the weights $\lambda=(\lambda(j))_\zs{1\le j\le n}$ belong to some finite set $\Lambda$ from $[0,1]^n$,
 $\wh{\theta}_\zs{j,\varepsilon}$ is defined in \eqref{sec:In.7} and $\phi_\zs{j}$ in \eqref{sec:In.5}. 
 Now we set
 \begin{equation}\label{sec:Mo.2}
\iota
=\card(\Lambda)
\quad\mbox{and}\quad
\vert\Lambda\vert_\zs{*}= \max_\zs{\lambda\in\Lambda}\,\sum^{n}_\zs{j=1}\,\Chi_\zs{\{\lambda_\zs{j}>0\}}
\,,
\end{equation}
where $\card(\Lambda)$ is the number of the vectors in $\Lambda$. In the sequel we assume that
 $\iota$ is a function of $\varepsilon>0$, i.e. $\iota=\iota(\varepsilon)$, such that for any $\b>0$
\begin{equation}
\label{cond-card-Lmbd-11}
\lim_\zs{\varepsilon\to 0}\,
\varepsilon^{\b} \iota(\varepsilon)=0\,.
\end{equation}
Now we chose the truncating parameter $\overline{a}_\zs{\varepsilon}$ as

 \begin{equation}\label{sec:Mo.2++}
\overline{a}_\zs{\varepsilon}=
\frac{\varepsilon}{\vert\Lambda\vert_\zs{*}}
\,.
\end{equation}

\noindent 
To choose a weight sequence $\lambda$ in the set $\Lambda$ we use the empirical  quadratic risk, defined as
$$
\Er_\varepsilon(\lambda) = \parallel \wh{S}_\lambda-S\parallel^2,
$$
which in our case is equal to
\begin{equation}\label{sec:Mo.3}
\Er_\varepsilon(\lambda) = \sum_\zs{j=1}^{n} \lambda^2(j) \wh{\theta}^2_\zs{j,\varepsilon} -2 \sum_\zs{j=1}^{n} \lambda(j) \wh{\theta}_\zs{j,\varepsilon}\theta_\zs{j}+ \sum_\zs{j=1}^{\infty} \theta^2_\zs{j}.
\end{equation}
Since the Fourier coefficients $(\theta_\zs{j})_\zs{j\ge\,1}$ are unknown, we replace
the terms $\wh{\theta}_\zs{j,\varepsilon}\theta_\zs{j}$ by  
\begin{equation}\label{sec:Mo.4}
\wt{\theta}_\zs{j,\varepsilon} = \wh{\theta}^2_\zs{j,\varepsilon} - \varepsilon^2 \wh{\varkappa}_\zs{\varepsilon}\,,
\end{equation}
where $\wh{\varkappa}_\zs{\varepsilon}$ is a some estimate for  the variance  parameter $\check{\varkappa}_\zs{Q}$ from \eqref{sec:In.3}. 
If it is known we set $\wh{\varkappa}_\zs{\varepsilon}=\check{\varkappa}_\zs{Q}$ if not this estimator will be prescribed  later.

\begin{remark}\label{Re.sec:Ma.00-1-0}
To understand the estimate  \eqref{sec:Mo.4} 
note that the natural way is to remplace in the production 
$\wh{\theta}_\zs{j}\theta_\zs{j}$ the unknown coefficient $\theta_\zs{j}$
 with its estimator $\wh{\theta}_\zs{j}$, so we obtain $\wh{\theta}^{2}_\zs{j}$. But this is not good estimator for the production
since in vue of \eqref{sec:In.8} we obtain $\E_\zs{Q}\, \wh{\theta}_\zs{j}\theta_\zs{j}=\theta^{2}_\zs{j}$, but  
$\E_\zs{Q}\, \wh{\theta}^{2}_\zs{j}=\theta^{2}_\zs{j}+\varepsilon^{2} \E_\zs{Q}\overline{\xi}^{2}_\zs{j}$. Therefore, 
to obtain unbiased estimator for the production  $\wh{\theta}_\zs{j}\theta_\zs{j}$ for $j\ge 2$ one needs to subtract the variance 
$\varepsilon^{2} \check{\varkappa}_\zs{Q}$ if $\check{\varkappa}_\zs{Q}$ is known and its estimate if non. This gives the form \eqref{sec:Mo.4}.
It should be noted also that we don't take into account the first term, i.e. the case $j=1$. But only one term has not sufficient influence in the total sum, i.e. it is negligible
in the empiric risk \eqref{sec:Mo.3}.
\end{remark}

\noindent 
Finally, to choose the weights we will minimize the following cost function
\begin{equation}\label{sec:Mo.5}
J_\varepsilon(\lambda)=\sum_\zs{j=1}^{n} \lambda^2(j) \wh{\theta}^2_\zs{j,\varepsilon} -2 \sum_\zs{j=1}^{n} \lambda(j)\wt{\theta}_\zs{j,\varepsilon} + 
\delta\,\wh{P}_\zs{\varepsilon}(\lambda)
\,,
\end{equation}
where $\delta>0$ is some  threshold which will be specified later and the penalty term
\begin{equation}\label{sec:Mo.6}
\wh{P}_\zs{\varepsilon}(\lambda)= \varepsilon^2  \wh{\varkappa}_\zs{\varepsilon} |\lambda|^{2}
\quad\mbox{and}\quad
|\lambda|^{2}=\sum^{n}_\zs{j=1}\,\lambda^{2}_\zs{j}\,.
\end{equation}
Note that, if the $\check{\varkappa}_\zs{Q}$ is known then the penalty term is defined as
\begin{equation}\label{sec:Mo.6+1}
P_\zs{\varepsilon}(\lambda)= \varepsilon^2  \,\check{\varkappa}_\zs{Q} |\lambda|^{2}
\,.
\end{equation}

\noindent 
We define the model selection procedure as
\begin{equation}\label{sec:Mo.9}
\wh{S}_* = \wh{S}_\zs{\hat \lambda}
\quad\mbox{and}\quad
\wh{\lambda}= \mbox{argmin}_\zs{\lambda\in\Lambda} J_\zs{\varepsilon}(\lambda)\,.
\end{equation}
We recall that the set $\Lambda$ is finite so $\hat \lambda$ exists. In the case when $\hat \lambda$ is not unique we take one of them.

Now we estimate the  variance parameter $\check{\varkappa}_\zs{Q}$
defined in  \eqref{sec:In.3}. To this end  for any $0<\varepsilon\le 1/\sqrt{3}$,
we set
\begin{equation}\label{sec:Mo.4-1-31-3}
\wh{\varkappa}_\zs{\varepsilon}=
\sum^n_\zs{j=[1/\varepsilon]+1}\,\wh{\tau}^2_\zs{j,\varepsilon}\,,
\quad n=[1/\varepsilon^{2}]\,,
\end{equation}
where $\wh{\tau}_\zs{j,\varepsilon}$ are the estimators for the Fourrier coefficients $\tau_\zs{j}$  with respect to the trigonometric basis \eqref{sec:In.5}, i.e.
\begin{equation}\label{sec:Mo.4-2-31-3}
\wh{\tau}_\zs{j,\varepsilon}=\,
\int^{1}_\zs{0}\,\Tr_\zs{j}(t)\d \check{y}_\zs{t}
\quad\mbox{and}\quad
\tau_\zs{j}=\,
\int^{1}_\zs{0}\,S(t)\Tr_\zs{j}(t)\d t
\,.
\end{equation}

\noindent 
We  study this estimator.

\begin{proposition}
\label{Pr.sec: MaPr.3}
Assume that in the model  \eqref{sec:In.1} the unknown function $S(\cdot)$ is continuously
differentiable.  Then, for any $0<\varepsilon\le 1/\sqrt{3}$
\begin{equation}
\label{sec: MaPr.5--1}
\E_\zs{Q}\,\vert \wh{\varkappa}_\zs{\varepsilon}-\check{\varkappa}_\zs{Q}\vert
\le\varepsilon
\Upsilon_\zs{Q}(S)
+
 \frac{\sqrt{6\check{\varkappa}_\zs{Q}}}{\vert\Lambda\vert_\zs{*}}
\,,
\end{equation}
where
$
\Upsilon_\zs{Q}(S)=4(\Vert\dot{S}\Vert+1)^{2}
\left(1+\sqrt{\check{\varkappa}_\zs{Q}}+2\check{\varkappa}_\zs{Q}+\sqrt{U_\zs{Q}}\right)$ and $\dot{S}$ is the derivative of the function $S$.
\end{proposition}
\noindent
The proof of this proposition is given in Section \ref{sec:Pr}. It is clear that in the case when $\vert\Lambda\vert_\zs{*}\le 1/\varepsilon$
we obtain that

\begin{equation}
\label{sec: MaPr.5--1+}
\E_\zs{Q}\,\vert \wh{\varkappa}_\zs{\varepsilon}
-
\check{\varkappa}_\zs{Q}
\vert
\le
 \frac{\Upsilon_\zs{Q}(S)+\sqrt{6\check{\varkappa}_\zs{Q}}}{\vert\Lambda\vert_\zs{*}}
\,.
\end{equation}

\begin{remark}\label{Re.sec:Ma.001++22--}
It should be noted that to estimate the parameter $\check{\varkappa}_\zs{Q}$
we use the equality \eqref{sec:In.8} for the Fourier coefficients $(\tau_\zs{j})_\zs{j\ge 1}$ with respect to the trigonometric basis
\eqref{sec:In.5}.  Moreover, as is shown in
Lemma A.6 in \cite{KonevPergamenshchikov2009a}
for any continuously differentiable function $S$
and 
for any $m\ge 1$
the sum $\sum_\zs{j\ge m}\,\tau^{2}_\zs{j}$ can be estimated from above in the explicite form. So,
taking this into account and  
 the properties  \eqref{sec: MaPr.2}
we the upper bound \eqref{sec: MaPr.5--1}.
\end{remark}

\noindent 
Now, we specify the weight coefficients
$(\lambda(j))_\zs{1\le j\le n}$. Consider a numerical grid of the form
\begin{equation}\label{sec:Ga.0}
\cA=\{1,\ldots,k^*\}\times\{r_1,\ldots,r_\zs{\m}\}\,,
\end{equation}
where  
$r_\zs{i}=i\,\varpi$ and
$\m=[1/\varpi^{2}]$.  We assume that both
the parameters $k^*\ge 1$ and $0<\varpi<1$ are functions of $\ve$, i.e.
$k^*=k^{*}_\zs{\ve}$ and $\varpi=\varpi_\zs{\ve}$ such that
\begin{equation}\label{sec:Ga.1}
\lim_\zs{\ve\to 0}\,
\left(
\dfrac{1}{k^{*}_\zs{\ve}}+
\dfrac{k^{*}_\zs{\ve}}{\vert\ln \ve\vert}
\right)
=0
\quad\mbox{and}
\quad 
\lim_\zs{\ve\to 0}\,
\left(\varpi_\zs{\ve}
+
\frac{\ve^{\b}}{\varpi_\zs{\ve}}
\right)
=0
\end{equation}
for any $\b>0$. One can take, for example, for $0<\ve<1$
\begin{equation}\label{sec:Ga.1-00}
\varpi_\zs{\ve}=\vert\ln\ve\vert^{-1}
\quad\mbox{and}\quad
k^{*}_\zs{\ve}=k^{*}_\zs{0}+\sqrt{\vert\ln \ve\vert}\,,
\end{equation}
where $k^{*}_\zs{0}\ge 0$ is some fixed constant.
 For each $\alpha=(\beta,r)\in\cA$, we introduce the weights $\lambda_\zs{\alpha}=(\lambda_\zs{\alpha}(j))_\zs{1\le j\le n}$ from $\bbr^{n}$
as
\begin{equation}\label{sec:Ga.2}
\lambda_\zs{\alpha}(j)=\Chi_\zs{\{1\le j<j_\zs{*}\}}+
\left(1-(j/\omega_\alpha)^\beta\right)\,
\Chi_\zs{\{ j_\zs{*}\le j\le \omega_\zs{\alpha}\}}
\,,
\end{equation}
where
$j_\zs{*}=j_\zs{*}(\alpha)=\left[\omega_\zs{\alpha}/\vert\ln\ve\vert\right]$,
$\omega_\zs{\alpha}=\d_\zs{\beta}\,(r\,\upsilon_\zs{\ve})^{1/(2\beta+1)}$, 
\begin{equation}\label{sec:Ga.2++tau}
\d_\zs{\beta}=
\left(\frac{(\beta+1)(2\beta+1)}{\pi^{2\beta}\beta}
\right)^{1/(2\beta+1)}
\,,
\qquad
\upsilon_\zs{\ve}=
\frac{1}{\ve^{2}\,\varsigma^{*}_\zs{\varepsilon}}
\end{equation}
and the threshold $\varsigma^{*}_\zs{\varepsilon}$ is introduced in \eqref{sec:Ex.01-1}.
Now we define the set $\Lambda$ 
as
\begin{equation}\label{sec:Ga.3}
\Lambda\,=\,\{\lambda_\zs{\alpha}\,,\,\alpha\in\cA\}\,.
\end{equation}

\noindent 
Note that in this case $\iota=k^* \m$ and the conditions \eqref{sec:Ga.1}
imply directly the property \eqref{cond-card-Lmbd-11}. Moreover, from
\eqref{sec:Ga.2}  we find that for any $\alpha\in\cA$
$$
\sum^{n}_\zs{j=1}\,\lambda_\zs{\alpha}(j)\le \omega_\zs{\alpha}
\le d_\zs{*}\,r^{1/3}_\zs{\m}\,\upsilon^{1/3}_\zs{\ve}
\quad\mbox{and}\quad
d_\zs{*}=\sup_\zs{\beta\ge 1}
\, \d_\zs{\beta}
\,.
$$
Therefore, the conditions
\eqref{sec:Ga.1} imply that for any $\b>0$

\begin{equation}\label{lim_Lambda-*}
 \lim_\zs{\varepsilon\to 0}\varepsilon^{2/3+\b} \vert\Lambda\vert_\zs{*}=0\,.
\end{equation}

\begin{remark}\label{Re.sec:Ma.2++22}
The parameters $\beta$ and $r$
are defined by the regularity of the unknown function $S$ 
(see for the details Remark \ref{Re.sec.Mrs.1} below).
It should be emphasized that
the weight coefficients defined by the set \eqref{sec:Ga.3} are
  used 
  by Konev and Pergamenshchikov
  in \cite{KonevPergamenshchikov2012, KonevPergamenshchikov2015}
 for continuous time regression
models  to show the asymptotic efficiency.
\end{remark}

\bigskip

\section{Main results}\label{sec:Mrs}

\subsection{Oracle inequalities}

First we set the following constant which will be used 
to describe the rest term 
in the oracle inequalities. We set

\begin{equation}\label{sec:OIn.1}
\Psi_\zs{Q,\varepsilon}=
(1+(\phi^{*})^{4})
\left( 
1+\varkappa^{2}_\zs{Q}
+\frac{1}{\check{\varkappa}_\zs{Q}}
\right)
\,\iota\,.
\end{equation}

\noindent 
We start with the sharp oracle inequalities.

\begin{theorem}\label{Th.sec:Mrs.1} 
Assume that for the model \eqref{sec:In.1} the condition
\eqref{sec:Ex.1-00}
 holds.
Then there exists a  constant $\l_\zs{*}>0$
such that
 for any $\varepsilon>0$ and $0 <\delta< 1/6$,  the estimator of $S$ given in \eqref{sec:Mo.9}
with the truncation parameter \eqref{sec:Mo.2++}
 satisfies the following oracle inequality
\begin{equation}\label{sec:OI.1}
\mathcal{R}_\zs{Q} (\wh{S}_*,S)\leq\frac{1+3\delta}{1-3\delta} \min_\zs{\lambda\in\Lambda} 
\mathcal{R}_\zs{Q} (\wh{S}_\lambda,S)+ \varepsilon^2\,\l_\zs{*}\,
 \frac{\Psi_\zs{Q,\varepsilon}
 + \vert\Lambda\vert_\zs{*}\, \E_\zs{S} |  \wh{\varkappa}_\zs{\varepsilon} -\check{\varkappa}_\zs{Q} |}{\delta}\,.
\end{equation}
\end{theorem}
\noindent
If the parameter $\check{\varkappa}_\zs{Q}$ is known,
we can simplify this inequality.
\begin{corollary}\label{Co.sec:OI.1} 
Assume that for the model \eqref{sec:In.1} the condition
\eqref{sec:Ex.1-00}
 holds.
If the variance parameter $\check{\varkappa}_\zs{Q}$ is known,
then there exists a constant $\l_\zs{*}>0$
such that
 for any $\varepsilon>0$ and  $0 <\delta< 1/6$,  the estimator of $S$ given in \eqref{sec:Mo.9}
with the truncation parameter \eqref{sec:Mo.2++}
 satisfies the following oracle inequality
\begin{equation}\label{sec:OIn.2}
\mathcal{R}_\zs{Q} (\wh{S}_*,S)\leq\frac{1+3\delta}{1-3\delta} \min_\zs{\lambda\in\Lambda} 
\mathcal{R}_\zs{Q} (\wh{S}_\lambda,S)+ \varepsilon^2\,\l_\zs{*}\,
 \frac{\Psi_\zs{Q,\varepsilon}}{\delta}\,.
\end{equation}
\end{corollary}

\begin{remark}\label{Re.sec:MainRes.1--00}
It should be noted that in the classical "signal+white noise" model,
i.e. when in the process \eqref{sec:In.1+1} the parameter $\varrho_\zs{1}=1$ and 
  the Levy measure $\Pi=0$,  we obtain $\check{\varkappa}_\zs{Q}=1$. Therefore,
 we can use the inequality \eqref{sec:OIn.2}.
\end{remark}

\noindent 
Using Proposition \ref{Pr.sec: MaPr.3} we can obtain the following inequality.

\begin{theorem}\label{Th.sec:Mrs.2--0}
Assume that for the model \eqref{sec:In.1} the condition
\eqref{sec:Ex.1-00}
 holds and the unknown signal $S(\cdot)$ is continuously
differentiable $[0,1]\to\bbr$ function. Then there exists a some constant $\l_\zs{*}>0$
such that for any
   $0 <\delta< 1/6$ and for any $\varepsilon>0$, for which $\vert\Lambda\vert_\zs{*}\le 1/\varepsilon$, 
    the estimator of $S$ given in \eqref{sec:Mo.9}
with the truncation parameter \eqref{sec:Mo.2++}
 satisfies the following oracle inequality
\begin{equation}\label{sec:Mrs.1+}
\mathcal{R}_\zs{Q} (\wh{S}_*,S)\leq\frac{1+3\delta}{1-3\delta} \min_\zs{\lambda\in\Lambda} 
\mathcal{R}_\zs{Q} (\wh{S}_\lambda,S)
+ \varepsilon^2\,\l_\zs{*}\,
 \frac{\Psi_\zs{Q,\varepsilon}(\Vert\dot{S}\Vert+1)^{2}}{\delta}\,.
\end{equation}
\end{theorem}

\bigskip

\noindent 
Now we study the robust risks defined in
\eqref{sec:robust-risks}
 for the procedure \eqref{sec:Mo.9}.

Moreover, we assume also that the upper bound for the basis functions  in \eqref{sec:In.3-00} 
may be dependent on $\varepsilon>0$, i.e. $\phi_\zs{*}=\phi_\zs{*}(\varepsilon)$, such that 
for any  $\b>0$
\begin{equation}\label{sec:Mrs.5-2}
\lim_\zs{n\to\infty}\,\varepsilon^{\b}\,
\phi_\zs{*}(\varepsilon)
=0\,.
\end{equation}

\begin{theorem}\label{Th.sec:Mrs.2}
Assume that for the model \eqref{sec:In.1} the condition
\eqref{sec:Ex.1-00}
 holds and the unknown function $S(\cdot)$ is continuously
differentiable. Then
for any  $0 <\delta< 1/6$
and 
for any $\varepsilon>0$ for which
 $\vert\Lambda\vert_\zs{*}\le 1/\varepsilon$, 
the
 robust risks for  the procedure \eqref{sec:Mo.9} 
with the truncation parameter \eqref{sec:Mo.2++}
 satisfies the following oracle inequality
\begin{equation}\label{sec:Mrs.6-25.3}
\cR^{*}_\zs{\varepsilon}(\wh{S}_*,S)\leq\frac{1+3\delta}{1-3\delta} \min_\zs{\lambda\in\Lambda}
 \cR^{*}_\zs{\varepsilon}(\wh{S}_\lambda,S)+\ve^{2}\,
\frac{\U^{*}_\zs{\ve}(S)}{\delta}
\,,
\end{equation}
where the term $\U^{*}_\zs{\varepsilon}(S)>0$ is such that under the conditions  \eqref{cond-card-Lmbd-11} and \eqref{sec:Mrs.5-2} 
for any $r>0$ and $\b>0$
\begin{equation}\label{sec:Mrs.7-25.3}
\lim_\zs{\ve\to 0}\,\ve^{\b}\,
\sup_\zs{\|\dot{S}\| \le r}
\,\U^{*}_\zs{\ve}(S)\,
=0
\,.
\end{equation}
\end{theorem}
\noindent Now taking into account the property
\eqref{lim_Lambda-*} we can deduce the following theorem for the procedure \eqref{sec:Mo.9}
with the weight coefficients \eqref{sec:Ga.3}.
\begin{theorem}\label{Th.sec:Mrs.++}
Assume that for the model \eqref{sec:In.1} the condition
\eqref{sec:Ex.1-00} holds.
Then the model selection procedure \eqref{sec:Mo.9}
constructed 
on the basis functions satisfying the condition \eqref{sec:Mrs.5-2}
and
through
 the weight coefficients \eqref{sec:Ga.3}
with the conditions
\eqref{sec:Ga.1}
 satisfies the oracle inequality
 \eqref{sec:Mrs.6-25.3}
 with the property
\eqref{sec:Mrs.7-25.3}. 
\end{theorem}

\begin{remark}\label{Re.sec:Mo.1}
Note that the similar sharp oracle inequalities were obtained  in the papers
\cite{GaltchoukPergamenshchikov2009a}
and
\cite{KonevPergamenshchikov2012}
for the model selection procedures based on the trigonometric basis functions \eqref{sec:In.5}.
In this paper we obtain these inequalities for the model selection procedures based
on any arbitrary orthogonal basic function in $\L_\zs{2}[0,1]$. We use the trigonometric functions 
only to estimate the noise parameter $\check{\varkappa}_\zs{Q}$.
\end{remark}

\bigskip

\subsection{Adaptive robust efficiency}\label{sec:AdpEst}

\noindent 
Now we study the asymptotically efficiency properties for the procedure
\eqref{sec:Mo.9}, \eqref{sec:Ga.2}
  with respect to the robust risks \eqref{sec:robust-risks} defined by the 
  distribution family \eqref{sec:Ex.01-1} -- \eqref{sec:Ex.01-2}.
   To this end we assume that the unknown function
 \eqref{sec:In.6}
 belongs to the following ellipsoid in $l_\zs{2}$
 \begin{equation}\label{sec:Ef.2}
W^{k}_\zs{\r}=\{S\in\,\L_\zs{2}[0,1]\,:\,
\sum_\zs{j=1}^{\infty}\,a_\zs{j}\,\theta^2_\zs{j}\,\le \r\}
 \end{equation}
where $a_\zs{j}=\sum^k_\zs{i=0}\left(2\pi [j/2]\right)^{2i}$.

 It is easy to see that in the case when the functions $(\phi_\zs{j})_\zs{j\ge 1}$ are trigonometric
 \eqref{sec:In.5},
 then this set coincides with the  Sobolev ball
\begin{equation}\label{sec:Ef.1}
W^{k}_\zs{\r}=\{f\in \,\C^{k}_\zs{per}[0,1]
\,:\,\sum_\zs{j=0}^k\,\|f^{(j)}\|^2\le \r\}\,,
 \end{equation}
where $\r>0$ and $ k\ge 1$ are some  parameters,
$\C^{k}_\zs{per}[0,1]$ is the set of
 $k$ times continuously differentiable functions
$f\,:\,[0,1]\to\bbr$ such that $f^{(i)}(0)=f^{(i)}(1)$ for all
$0\le i \le k$. Similarly to
\cite{KonevPergamenshchikov2012, KonevPergamenshchikov2015} 
we will  show here that the asymptotic sharp lower bound 
 for the robust risk \eqref{sec:robust-risks}
is given by
\begin{equation}\label{sec:Ef.3}
l_\zs{*}(\r)
=
\,
\left((2k+1)\r\right)^{1/(2k+1)}\,
\left(
\frac{k}{(k+1)\pi} \right)^{2k/(2k+1)}\,.
\end{equation}
\noindent Note that this is  the well-known Pinsker constant
obtained for the nonadaptive filtration problem in ``signal +
small white noise'' model
 (see, for example, \cite{Pinsker1981}). Let $\cS_\zs{\ve}$ be the set of all estimators $\wh{S}_\zs{\ve}$
measurable with respect to the sigma-algebra
$\sigma\{y_\zs{t}\,,\,0\le t\le 1\}$
 generated by the process \eqref{sec:In.1}.

\begin{theorem}\label{Th.sec:Ef.1} 
For the distribution family \eqref{sec:Ex.01-1} -- \eqref{sec:Ex.01-2}.
the robust risk admits the following lower bound
 \begin{equation}\label{sec:Ef.4}
\liminf_\zs{\ve\to 0}\,
\upsilon
^{2k/(2k+1)}_\zs{\ve}\,
 \inf_\zs{\wh{S}_\zs{\ve}\in\cS_\zs{\ve}}\,\,
\sup_\zs{S\in W^{k}_\zs{r}} \,\cR^{*}_\zs{\varepsilon}(\wh{S}_\zs{\varepsilon},S) \ge
\,l_\zs{*}(\r)\,,
 \end{equation}
where the rate $\upsilon_\zs{\ve}$ is given in \eqref{sec:Ga.2++tau}, i.e. $\upsilon_\zs{\ve}=
\left(\ve^{2}\,\varsigma^{*}_\zs{\varepsilon}\right)^{-1}$.
\end{theorem}

\noindent 
We set the parameter $\delta$ in \eqref{sec:Mo.5} as function of $\varepsilon$, i.e.
$\delta=\delta_\zs{\varepsilon}$ is such that
 \begin{equation}\label{sec:Ef.4-01}
\lim_\zs{\ve\to 0}\,\delta_\zs{\ve}=0
\quad\mbox{and}\quad
\lim_\zs{\ve\to 0}\,\varepsilon^{-\b}\,\delta_\zs{\varepsilon}=+\infty
 \end{equation}
for any $\b>0$. For example, we can take $\delta_\zs{\varepsilon}=(6+\vert\ln \varepsilon\vert)^{-1}$.

\begin{theorem}\label{Th.sec:Ef.2}
Assume that for the model \eqref{sec:In.1} the condition
\eqref{sec:Ex.1-00} holds.
Then the model selection procedure \eqref{sec:Mo.9}
constructed 
on the basis functions satisfying the condition \eqref{sec:Mrs.5-2}
and
through
 the weight coefficients \eqref{sec:Ga.3}
with the conditions
\eqref{sec:Ga.1}
 admits 
 the following asymptotic upper bound
 \begin{equation}\label{sec:Ef.5}
\limsup_\zs{\varepsilon\to 0}\,
\upsilon^{2k/(2k+1)}_\zs{\varepsilon}\,
 \sup_\zs{S\in W^k_r}\,
\cR^{*}_\zs{\varepsilon}(\wh{S}_\zs{*},S) \le  
l_\zs{*}(\r)
\,.
 \end{equation}
\end{theorem}

\medskip

\noindent Theorem~\ref{Th.sec:Ef.1} and Theorem~\ref{Th.sec:Ef.2} imply the following result
\begin{corollary}\label{Co.sec:Mr.1}
Under the conditions Theorem~\ref{Th.sec:Ef.2}
\begin{equation}\label{sec:Ef.6}
\lim_\zs{\varepsilon\to 0}\,
\upsilon
^{2k/(2k+1)}_\zs{\varepsilon}\,
 \inf_\zs{\wh{S}_\zs{\varepsilon}\in\cS_\zs{\varepsilon}}\,\,
\sup_\zs{S\in W^{k}_\zs{r}} \,\cR^{*}_\zs{\varepsilon}(\wh{S}_\zs{\varepsilon},S)
= l_\zs{*}(\r)\,.
 \end{equation}
\end{corollary}

 \begin{remark}
 \label{Re.sec.Mrs.1}
It should be noted 
(see, for example, Pinsker \cite{Pinsker1981}). that if the parameters $k$ and $r$ of the Sobolev ball
\eqref{sec:Ef.2}
 are known, then to obtain the efficient estimation it suffice to chose the weight least square estimator 
 \eqref{sec:Mo.1++1} with the weights \eqref{sec:Ga.2}
 and $\alpha=(k,r)$. In the adaptive estimation case, i.e. when these parameters are unknown we propose to use the
 selection model procedure for the family $(\wh{S}_\zs{\lambda})_\zs{\lambda\in\Lambda}$ which contains the efficient estimator.
Then, the efficiency property is provided through the sharp oracle inequalities.  
Moreover, note also that
 the optimal (minimax) risk convergence rate
for the Sobolev ball $W^{k}_\zs{r}$
 is $\varepsilon^{-4k/(2k+1)}$.
We see here that the efficient robust rate is
 $\upsilon^{2k/(2k+1)}_\zs{\varepsilon}$, i.e. if  the distribution upper bound  $\varsigma^{*}_\zs{\varepsilon}\to 0$  as $n\to\infty$
 we obtain the more rapid rate with respect to $\varepsilon^{-4k/(2k+1)}$, and if $\varsigma^{*}_\zs{\varepsilon}\to \infty$  as 
 $\varepsilon\to 0$
we obtain the more slow rate. In the case when $\varsigma^{*}_\zs{\varepsilon}$ is constant the robuste rate is the same as the classical non robuste convergence rate.
\end{remark}

\bigskip
\section{The van Trees inequality for L\'evy processes.}\label{sec: VanTreesIn__}

In this section we consider the following continuous time
 parametric regression model 
\begin{equation}\label{sec:App.5}
  \d y_t=S(t,\theta)\d t+\d \xi_t\,,
  \quad 0\le t\le 1\,,
 \end{equation}
where  $S(t,\theta)=\sum^{d}_\zs{i=1}\,\theta_\zs{i}\,\phi_\zs{i}(t)$
with the unknown parameters
$\theta=(\theta_\zs{1},\ldots,\theta_\zs{d})'$
and the process
$(\xi_\zs{t})_\zs{0\le t\le 1}$
is defined in
\eqref{sec:In.1+1}. 
Note now that according to Proposition \ref{Pr.sec:App.1++}
the distribution
$\P_\zs{\theta}$
of the process \eqref{sec:App.5} 
 is absolutely continuous with respect to the 
$\P_\zs{\xi}$ on $\D[0,1]$
and the corresponding Radon-Nikodym derivative is
\begin{equation}\label{sec:App.7}
f(x,\theta)=
\frac{\d\P_\zs{\theta}}{\d\P_\zs{\xi}}(x)=
\exp\left\{\int^{1}_\zs{0}\,\frac{S(t,\theta)}{\varrho^{2}_\zs{1}}\,\d x^{c}_\zs{t}
-\,\int^{1}_\zs{0}\,
\frac{S^{2}(t,\theta)}{2\varrho^{2}_\zs{1}}\,
\d t
\right\}
\,,
\end{equation}
where $(x^{c}_\zs{t})_\zs{0\le t\le T}$ is the continuous part of the process $(x_\zs{t})_\zs{0\le t\le T}$ in $\D[0,T]$, i.e.
$$
x^{c}_\zs{t}=
x_\zs{t}
-
\int^{t}_\zs{0}\,\int_\zs{\bbr}\,v\,\left( 
\mu_\zs{x}(\d s\,,\d v)
-
\Pi(\d v)\d s
\right)
$$
and for any $t>0$ and any measurable  $\Gamma$ from $\bbr\setminus \{0\}$
$$
\mu_\zs{x}([0,t],\Gamma)
=\sum_\zs{0\le s\le t}\,\Chi_\zs{\{\Delta x_\zs{s}\in\varrho_\zs{2}\Gamma\}}
\,.
$$


Let $\Phi$ be a prior density on $\bbr^d$ having
the following form:
$$
\Phi(\theta)=\Phi(\theta_1,\ldots,\theta_d)=\prod_{j=1}^d\varphi_\zs{j}(\theta_\zs{j})\,,
$$
where $\varphi_\zs{j}$ is some continuously differentiable density in $\bbr$. 
Moreover, let $g(\theta)$ be a continuously differentiable $\bbr^d\to \bbr$ function such that,
for each $1\le j\le d$,
\begin{equation}\label{sec:App.8}
\lim_\zs{|\theta_\zs{j}|\to\infty}\,
g(\theta)\,\varphi_\zs{j}(\theta_\zs{j})=0
\quad\mbox{and}\quad
\int_\zs{\bbr^d}\,|g^{\prime}_\zs{j}(\theta)|\,\Phi(\theta)\,\d \theta
<\infty\,,
\end{equation}
where
$$
g^{\prime}_\zs{j}(\theta)=\frac{\partial g(\theta)}{\partial\theta_\zs{j}}\,.
$$
For any $\cB(\cX)\times\cB(\bbr^d)-$
measurable integrable function $H=H(x,\theta)$ we denote
\begin{align*}
\wt{\E}\,H&=\int_{\bbr^d}\,
\int_\zs{\cX}\,H(x,\theta)\,\d \P_\zs{\theta}\,\Phi(\theta) \d \theta
\\[2mm]&
=
\int_{\bbr^d}\,\int_\zs{\cX}\,
H(x,\theta)\,f(x,\theta)\,\Phi(\theta)\d \P_\zs{\xi}(x)\, \d \theta\,,
\end{align*}
where $\cX=\D[0,1]$.

\begin{proposition}\label{Pr.sec:VanTreesIn++}
For any $\cF^y=\sigma\{y_\zs{t}\,0\le t\le 1\}$-measurable square integrable function $\wh{g}$
 and for any $1\le j\le d$, the following inequality holds
\begin{equation}
\label{sec:VanTrsInq_11_22++0}
\wt{\E}(\wh{g}-g(\theta))^2\ge
\frac{\Lambda^2_\zs{j}}{\Vert \phi_\zs{j}\Vert^{2}\varrho^{-2}_\zs{1}+I_\zs{j}}
\,,
\end{equation}
where
$$
\lambda_\zs{j}=\int_\zs{\bbr^d}\,g^{\prime}_\zs{j}(\theta)\,\Phi(\theta)\,\d \theta
\quad\mbox{and}\quad
I_\zs{j}=\int_\zs{\bbr}\,\frac{\dot{\varphi}^2_\zs{j}(z)}{\varphi_\zs{j}(z)}\,\d z\,.
$$
\end{proposition}
\noindent {\bf  Proof.}
First of all note that, the density \eqref{sec:App.7} 
on the process $\xi$
is bounded 
with respect to $\theta_\zs{j}\in\bbr$ and for any $1\le j\le d$
$$
\limsup_\zs{|\theta_\zs{j}|\to\infty}\,f(\xi,\theta)\,=\,0\,.
\quad\quad\mbox{a.s.}
$$
Now, we set
$$
\wt{\Phi}_\zs{j}=\wt{\Phi}_\zs{j}(x,\theta)=
\frac{\partial\,(f(x,\theta)\Phi(\theta))/\partial\theta_\zs{j}}{f(x,\theta)\Phi(\theta)}
 \,.
$$
Taking into account the condition \eqref{sec:App.8} and
integrating by parts yield 
\begin{align*}
\wt{\E}&\left(((\wh{g}-g(\theta))\wt{\Phi}_\zs{j}\right)
=\int_{\cX\times\bbr^d}\,((\wh{g}(x)-g(\theta))\frac{\partial}{\partial\theta_\zs{j}}
\left(f(x,\theta)\Phi(\theta)\right)\d \theta\,\P_\zs{\xi}(\d x)\\[2mm]
&=\int_{\cX\times\bbr^{d-1}}\left(\int_{\bbr}\,
g^{\prime}_\zs{j}(\theta)\,
f(x,\theta)\Phi(\theta)\d \theta_\zs{j}\right)\left(\prod_{i\neq j}\d \theta_i\right)\,\P_\zs{\xi}(\d x)
=\lambda_\zs{j}\,.
\end{align*}
Now by the Bouniakovskii-Cauchy-Schwarz inequality
we obtain the following lower bound for the quadratic risk
$$
\wt{\E}((\wh{g}-g(\theta))^2\ge
\frac{\Lambda^2_\zs{j}}{\wt{\E}\wt{\Phi}_\zs{j}^2}\,.
$$
To study the denominator in the left hand of this 
inequality note that in view of the reprentation
 \eqref{sec:App.7}  
$$
\frac{1}{f(y,\theta)}
\frac{\partial\,f(y,\theta)}{\partial\theta_\zs{j}}
=\frac{1}{\varrho_\zs{1}}\,
\int^{1}_\zs{0}\,\phi_\zs{j}(t)\,\d w_\zs{t}\,.
$$
Therefore, for each $\theta\in\bbr^d$,
$$
\E_\zs{\theta}\,
\frac{1}{f(y,\theta)}
\frac{\partial\,f(y,\theta)}{\partial\theta_\zs{j}}
\,
=0
$$
and
$$
\E_\zs{\theta}\,
\left(
\frac{1}{f(y,\theta)}
\frac{\partial\,f(y,\theta)}{\partial\theta_\zs{j}}
\right)^2
=\,
\frac{1}{\varrho^{2}_\zs{1}}
\int^{1}_\zs{0}\,\phi^2_\zs{j}(t)\d t
=
\frac{1}{\varrho^{2}_\zs{1}}
\Vert\phi_\zs{j}\Vert^{2}
\,.
$$
Taking into account that
$$
\wt{\Phi}_\zs{j}=
\frac{1}{f(x,\theta)}
\frac{\partial\,f(x,\theta)}{\partial\theta_\zs{j}}
+
\frac{1}{\Phi(\theta)}
\frac{\partial\,\Phi(\theta))}{\partial\theta_\zs{j}}
 \,,
$$
we get
$$
\wt{\E} \wt{\Phi}_\zs{j}^2=
\frac{1}{\varrho^{2}_\zs{1}}\,\Vert\phi_\zs{j}\Vert^{2}
+\,I_\zs{j}\,.
$$
Hence
Proposition~\ref{Pr.sec:VanTreesIn++}.
\endproof

 \begin{remark}
 \label{Re.sec.VanTreesIn++}
Note that, the lower bound \eqref{sec:VanTrsInq_11_22++0}
is an extension for the van Trees inequality used
for the "signal+white noise" model (see, for example, the inequality (A.5) in \cite{KonevPergamenshchikov2009b}). 
\end{remark}

\bigskip

\section{Lower bound}\label{sec:Lobn}

Firstly, note, that for any fixed $Q\in \cQ^{*}_\zs{\varepsilon}$
\begin{equation}\label{sec:Lo.1-0}
\sup_\zs{S\in W^{k}_\zs{r}}\,\cR^{*}_\zs{\varepsilon}(\wh{S}_\zs{\varepsilon},S)
\ge 
\,
\sup_\zs{S\in W^{k}_\zs{r}}\,
\cR_\zs{Q}(\wh{S}_\zs{\varepsilon},S)
\,.
\end{equation}
Now for any fixed $0<\check{\gamma}<1$ we set 
\begin{equation}\label{sec:Lo.1}
d=d_\zs{\varepsilon}=\left[\frac{k+1}{k}\upsilon^{1/(2k+1)}_\zs{\varepsilon}\,l_\zs{*}(r_\zs{0})\right]
\quad\mbox{and}\quad
\r_\zs{0}=(1-\check{\gamma})\r\,.
\end{equation}
Using this definition  
we introduce the parametric family $(S_\zs{z})_\zs{z\in\bbr^{d}}$ 
as
\begin{equation}\label{sec:Lo.4}
S_\zs{z}(x)=\sum_{j=1}^{d}\,z_\zs{j}\,\phi_\zs{j}(x)
\,.
\end{equation}

\noindent 
To define the bayesian risk we  choose a prior distribution on $\bbr^{d}$ 
as
\begin{equation}\label{sec:Lo.5}
\kappa=(\kappa_\zs{j})_\zs{1\le j\le d}
\quad\mbox{and}\quad
\kappa_\zs{j}=s_\zs{j}\,\eta_\zs{j}\,,
\end{equation}
where $\eta_\zs{j}$ are i.i.d. gaussian $\cN(0,1)$ random variables and
the coefficients 
$$
s_\zs{j}=\sqrt{\frac{s^*_\zs{j}}{v_\zs{\varepsilon}}}
\quad\mbox{and}\quad
s^{*}_\zs{j}\,
=
\left( \frac{d}{j}
\right)^{k}
-
1
\,.
$$
Denoting by $\mu_\zs{\kappa}$ the distribution of the random variables 
$(\kappa_\zs{j})_\zs{1\le j\le d}$ on $\bbr^{d}$
  we introduce
 the  Bayes risk as
\begin{equation}\label{sec:Lo.5_11-1}
\wt{\cR}_\zs{Q}(\wh{S})=
\int_\zs{\bbr^d}\cR_\zs{Q}(\wh{S},S_\zs{z})\,
\mu_\zs{\kappa}(\d z)\,.
\end{equation}
Furthermore, for any function $f\in\L_\zs{2}[0,1]$, we denote by $\p(f)$ its projection
in $\L_\zs{2}[0,1]$
 onto 
 $W_\zs{k,r}$, i.e.
 $$
\Vert f-\p(f)\Vert=\inf_\zs{h\in W^{k}_\zs{r}}\,\Vert f-h \Vert 
\,.
$$
Since $W^{k}_\zs{r}$ is a convex and closed set in $\L_\zs{2}[0,1]$, 
this projector exists and is unique for any function $f\in\L_\zs{2}[0,1]$ and,  moreover, 
$$
\|f-h\|^2\ge\|\p(f)-h\|^2
\quad\mbox{for any}\quad 
h\in W^{k}_\zs{r}
\,.
$$
So, setting $\wh{\p}=\p(\wh{S})$, we obtain that
$$
\sup_\zs{S\in W^{k}_\zs{r}}\,\cR(\wh{S},S)
\ge\,
\int_\zs{\{z\in\bbr^d\,:\,S_\zs{z}\in W^{k}_\zs{r}\}}\,
\E_\zs{S_\zs{z}}\|\wh{\p}-S_\zs{z}\|^2\,\mu_{\kappa}(\d z)
\,.
$$
Taking into account now that $\|\wh{\p}\|^2\le \r$ 
 we obtain 
\begin{equation}\label{sec:Lo.12}
\sup_\zs{S\in W^{k}_\zs{r}}\,
\cR_\zs{Q}(\wh{S},S)
\,
\ge\,
\wt{\cR}_\zs{Q}(\wh{\p})
-2\,
\Delta_\zs{\varepsilon}
\end{equation}
and
$$
\Delta_\zs{\varepsilon}=
\int_\zs{ \{z\in\bbr^d\,:\,S_\zs{z}\notin W_\zs{k,\r}\}}\,
\,
(\r+\|S_\zs{z}\|^2)\,
\mu_\zs{\kappa}(\d z)
\,.
$$
Therefore, in view of \eqref{sec:Lo.1-0},
\begin{equation}\label{sec:Lo.1+10}
\sup_\zs{S\in W_\zs{k,\r}}\,\cR^{*}_\zs{\varepsilon}(\wh{S}_\zs{\varepsilon},S)
\ge 
\,
\sup_\zs{Q\in\cQ^{*}_\zs{\varepsilon}}\,
\wt{\cR}_\zs{Q}(\wh{\p})
-2\,
\Delta_\zs{\varepsilon}
\,.
\end{equation}
As to the last term in this inequality, in Appendix we show that
for any $\b>0$
\begin{equation}\label{sec:Lo.12_11.1}
\lim_\zs{\varepsilon\to 0}\,
\varepsilon^{-\b}\,
\Delta_\zs{\varepsilon}=0\,.
\end{equation}
Now it is easy to see that
$$
\|\wh{\p}-S_z\|^2 \ge
\sum_{j=1}^{d}\,
(\wh{z}_\zs{j}-z_\zs{j})^2
\,,
$$
where $\wh{z}_\zs{j}=\int^{1}_\zs{0}\,\wh{\p}(t)\,\phi_\zs{j}(t)\d t.$ So, in view of Proposition~\ref{Pr.sec:VanTreesIn++}
and reminding
that 
$\upsilon_\zs{\varepsilon}=\varepsilon^{-2}/\varsigma^{*}_\zs{\varepsilon}$
 we obtain 
\begin{align*}
\sup_\zs{Q\in\cQ^{*}_\zs{\varepsilon}}
\wt{\cR}_\zs{Q}(\wh{\p})\,
&\ge\,
\sup_\zs{0<\varrho^{2}_\zs{1}\le \varsigma^{*}_\zs{\varepsilon}}
\sum_{j=1}^{d}\,\frac{1}
{\varepsilon^{-2}\,\varrho^{-2}_\zs{1}+v_\zs{\varepsilon}\,(s^{*}_\zs{j})^{-1}}
\\[2mm]
&
=
\frac{1}{v_\zs{\varepsilon}}\,
\sum_{j=1}^{d}\,\frac{s^{*}_\zs{j}}
{s^{*}_\zs{j}+\,1}
=
\frac{1}{v_\zs{\varepsilon}}\,
\sum_{j=1}^{d}\,
\left(
1
-
\frac{j^k}{d^k_\zs{\varepsilon}}
\right)
\,.
\end{align*}
Therefore, using now the definition \eqref{sec:Lo.1},  the inequality
\eqref{sec:Lo.1+10} and the limit
\eqref{sec:Lo.12_11.1},
we obtain that
$$
\liminf_\zs{n\to\infty}\inf_\zs{\wh{S}\in\Pi_\zs{\varepsilon}}\,v^{\frac{2k}{2k+1}}_\zs{\varepsilon}\,
\sup_\zs{S\in W_\zs{k,r}}\,\cR^{*}_\zs{\varepsilon}(\wh{S}_\zs{\varepsilon},S)
\ge\,
(1-\check{\gamma})^{\frac{1}{2k+1}}\,
l_\zs{*}(\r)\,.
$$
Taking here limit as $\check{\gamma}\to 0$ implies Theorem~\ref{Th.sec:Ef.1} . 
\endproof

\section{Upper bound}\label{sec:UpBn}

First of all
 we recall  the Novikov  inequalities, \cite{Novikov1975}, also referred to as the Bichteler--Jacod inequalities, see  
\cite{BichtelerJacod1983, MarinelliRockner2014},  providing bounds of the moments of the supremum of purely discontinuous local martingales for $p\ge 2$

\begin{equation}
\label{Novikov++}
\E\sup_{t\le 1}|g*(\mu-\wt{\mu})_\zs{t}|^{p}\le C^{*}_\zs{p}
 \left( 
 \E\,\big (|g|^{2}*\wt{\mu}_\zs{1}\big)^{p/2}
 +
 \E\,\big (|g|^{p}*\wt{\mu}_\zs{1}\big)
\right)\,,
\end{equation}
where $C^{*}_\zs{p}$ is some positive constant.

\subsection{Known smoothness}

First we suppose that the parameters $k\ge 1$, $\r>0$ in \eqref{sec:Ef.1} and
$\varsigma^{*}_\zs{\varepsilon}$ in \eqref{sec:Ex.01-1}  are known. Let the family of admissible
weighted least squares estimates
$(\wh{S}_\zs{\lambda})_\zs{\lambda\in\Lambda}$
 given
by
\eqref{sec:Ga.3}. Consider the pair
$$
\check{\alpha}=(k,\check{r})
\quad\mbox{and}\quad
\check{r}=\varpi \left[\r/\varpi \right]
\,.
$$
Denote the corresponding estimate as
\begin{equation}\label{sec:Up.1}
\check{S}=\wh{S}_\zs{\check{\lambda}}
\quad\mbox{and}\quad
\check{\lambda}=\lambda_\zs{\check{\alpha}}\,.
\end{equation}
\noindent Note that for sufficiently small $\varepsilon$ the pair
$\check{\alpha}$ belongs to the set \eqref{sec:Ga.0}.

\begin{theorem}\label{Th.sec:Up.1}
The estimator $\check{S}$
admits 
the following asymptotic upper
bound
\begin{equation}\label{Sec:Up.2}
\limsup_\zs{\varepsilon\to 0}\,
\upsilon^{2k/(2k+1)}_\zs{\varepsilon}
\,
\sup_{S\in W^{k}_\zs{\r}}\,
\cR^{*}_\zs{\varepsilon}\,(\check{S},S)\,
\le l_\zs{*}(\r)\,.
\end{equation}
\end{theorem}
\noindent {\bf Proof.}
Substituting \eqref{sec:In.8}
and taking into account the definition
\eqref{sec:Up.1}
one gets
$$
\|\check{S}-S\|^2
=\sum_{j=1}^{\infty}\,(1-\check{\lambda}(j))^2\,\theta^2_\zs{j}-2\check{M}_\zs{\varepsilon}
 +\varepsilon^{2}\,
 \sum_{j=1}^{\infty}\,\check{\lambda}^2(j)\,\check{\xi}^2_\zs{j}\,,
$$
where
$
\check{M}_\zs{\varepsilon}\,=\,\varepsilon\,
\sum_{j=1}^{\infty}\,(1\,-\,\check{\lambda}(j))\,\check{\lambda}(j)\,\theta_\zs{j}\,\overline{\xi}_\zs{j}$.
Note now that
for any $Q\in\cQ^{*}_\zs{\varepsilon}$
the expectation
$\E_\zs{Q,S}\,\check{M}_\zs{\varepsilon}=0$
and, in view of the upper bound \eqref{sec:Ex.01-1},
$$
\sup_\zs{Q\in\cQ^{*}_\zs{\varepsilon}}\,
\E_\zs{Q,S}
\sum_{j=1}^{\infty}\,\check{\lambda}^2(j)
\,\overline{\xi}^{2}_\zs{j}\le \varsigma^{*}_\zs{\varepsilon}\sum_{j=1}^{\infty}\,\check{\lambda}^2(j)\,.
$$
Therefore,
\begin{equation}\label{sec:Up.3}
\cR^{*}_\zs{\varepsilon}(\check{S},S)
\,\le\,
\sum_{j=\check{j}_\zs{*}}^{\infty}\,(1-\check{\lambda}(j))^2\,\theta^2_\zs{j}+
\frac{1}{\upsilon_\zs{\ve}}\sum_{j=1}^{\infty}\,\check{\lambda}^2(j)
\,,
\end{equation}
where $\check{j}_\zs{*}=j_\zs{*}(\check{\alpha})$. Setting
$$
\u_\zs{\varepsilon}= \upsilon^{2k/(2k+1)}_\zs{\ve} \sup_\zs{j\ge
\check{j}_\zs{*}}(1-\check{\lambda}(j))^2/a_\zs{j}\,,
$$
 we obtain that for each $S\in W^{k}_\zs{r}$
$$
\Upsilon_\zs{1,\varepsilon}(S)=\,\upsilon^{2k/(2k+1)}_\zs{\ve}\,
\sum_{j=\check{\iota}}^{\infty}\,(1-\check{\lambda}(j))^2\,\theta^2_\zs{j}
\,
\le\,\u_\zs{\varepsilon}\,
\sum_{j=\check{\iota}}^{\infty}\,a_\zs{j}\,
\theta^2_\zs{j}
\,
\le\,\u_\zs{\varepsilon}\,r\,.
$$
Taking into account that $\check{r}\to r$, we obtain that
$$
\limsup_\zs{\varepsilon\to 0}\,
\sup_\zs{S\in W^{k}_\zs{r}}\,
\Upsilon_\zs{1,\varepsilon}(S)
\,\le\,
\frac{r^{1/(2k+1)}}{\pi^{2k}(\d_\zs{k})^{2k/(2k+1)}}
:=\Upsilon^{*}_\zs{1}
\,.
$$
To estimate the last term in the right hand of \eqref{sec:Up.3}, we set
$$
\Upsilon_\zs{2,\varepsilon}=\,
\frac{1}{\upsilon^{1/(2k+1)}_\zs{\ve}}
\sum_{j=1}^{n}\,\check{\lambda}^2(j)\,.
$$
It is easy to check that
$$
\limsup_\zs{\varepsilon\to 0}\,
\Upsilon_\zs{2,\varepsilon}\le
\frac{2(r\d_\zs{k})^{1/(2k+1)}\,k^2}{(k+1)(2k+1)}
:=\Upsilon^{*}_\zs{2}
\,.
$$
Therefore, taking into account that by the definition of the Pinsker
constant in \eqref{sec:Ef.3}
$\Upsilon^{*}_\zs{1}
+
\Upsilon^{*}_\zs{2}
=
l_\zs{*}(\r)$,
we arrive at the inequality
$$
\lim_\zs{\varepsilon\to 0}\,
\upsilon^{2k/(2k+1)}_\zs{\varepsilon}
\sup_\zs{S\in W^{k}_\zs{r}}
\,\cR^{*}_\zs{\varepsilon}(\check{S},S)
\le\,l_\zs{*}(\r)\,.
$$
Hence Theorem~\ref{Th.sec:Up.1}.
\endproof

\subsection{Unknown smoothness}

Combining  Theorem~\ref{Th.sec:Mrs.++}
and
Theorem~\ref{Th.sec:Up.1} 
yields Theorem~\ref{Th.sec:Ef.2}.
\endproof

\bigskip

\section{Signals number  detection }\label{sec:NumSgn}

In this section we consider  the estimation problem for the signals  number
in the multi-path connection channel.
In the  framework of the statistical radio-physics models
we study
 the telecommunication system 
in which  
we observe 
the summarized signal
in the multi-path channel
with  noise on the time interval $[0,1]$:
$$
 y_\zs{t}= \sum^{q}_\zs{j=1}\theta_\zs{j}\phi_\zs{j}(t)
+
\nu_\zs{t}\,,
\quad 0\le t\le 1\,,
$$
where $(\nu_\zs{t})_\zs{t\ge 0}$  is the gaussian white noise.
The energetic parameters $(\theta_\zs{j})_\zs{j\ge 1}$,
and the number of  signals $q$
are unknown and
the signals $(\phi_\zs{j})_\zs{j\ge 1}$ are known orthonormal functions,
 i.e.
$$
\int^{1}_\zs{0}\,\phi_\zs{i}(t)\,\phi_\zs{j}(t)\,\d t=\Chi_\zs{\{i\neq j\}}
\,.
$$
 The problem is to estimate $q$ when
signal/noise ratio goes to infinity. To describe this problem in a mathematical framework  one has to use the following stochastic differential equation
\begin{equation}\label{sec:In.1-00_}
 \d y_\zs{t}= \left(\sum^{q}_\zs{j=1}\theta_\zs{j}\phi_\zs{j}(t)\right)\d t
+\varepsilon
\d w_\zs{t}
\,,
\end{equation}
where $(w_\zs{t})_\zs{t\ge 0}$ is the standard brownian motion and the parameter $\varepsilon>0$
is the noise intensity.
 We study this model when the signal/noise ration goes to infinity, i.e. $\varepsilon\to 0$.
 The logarithm of the likelihood ratio for  model  \eqref{sec:In.1-00_} can be represented as
$$
\ln L_\zs{\varepsilon}=
\frac{1}{\varepsilon^{2}}
 \sum^{q}_\zs{j=1}\theta_\zs{j}
\int^{1}_\zs{0}
\phi_\zs{j}(t)\d y_\zs{t}
-
\frac{1}{2\varepsilon^{2}}\sum^{q}_\zs{j=1}\,\theta^{2}_\zs{j}\,.
$$
If we  try to construct the maximum likelihood estimators for
  $(\theta_\zs{j})_\zs{1\le j\le q}$ and $q$, then we obtain that 
$$
\max_\zs{1\le q\le q_\zs{*}}
\max_\zs{\theta_\zs{j}}\ln L_\zs{\varepsilon}
=\frac{1}{2\varepsilon^{2}}
\sum^{q_\zs{*}}_\zs{j=1}\,\left(\int^{1}_\zs{0}
\phi_\zs{j}(t)\d y_\zs{t} 
\right)^{2}
\,.
$$
Therefore, the maximum likelihood estimation for $\wh{q}=q^{*}$. So, if $q^{*}=\infty$ we obtain that $\wh{q}=\infty$. 
Thus, this estimator gives nothing, i.e. it
 does not work. 
By these reasons  we propose to study the estimation problem for $q$ for the process \eqref{sec:In.1-00_} in a nonparametric setting
and to  apply  the model selection procedure \eqref{sec:Mo.9}. 
To this end
we consider the model
\eqref{sec:In.1}
 with the  unknown function $S$ defined as
\begin{equation}\label{NumSgn.1}
S(t)=\sum^{q}_\zs{j=1}\,\theta_\zs{j}\,\phi_\zs{j}(t)
\,.
\end{equation}
For this problem we use the LSE family $(\wh{S}_\zs{d})_\zs{1\le d\le m}$
defined as
\begin{equation}\label{sec:NumSgn.1}
\wh{S}_\zs{d} (x) = \sum_\zs{j=1}^{d} \, \wh{\theta}_\zs{j,\varepsilon} \phi_\zs{j}(x)\,.
\end{equation}
This estimate can be obtained  from \eqref{sec:Mo.1}
with the weights $\lambda_\zs{d}(j)=\chi{\{j\le d\}}$.
The number of estimators $\iota$ satisfies the condition
\eqref{cond-card-Lmbd-11}.
As a risk for the signals  number we use 
\begin{equation}\label{sec:NumSgn.3}
\D_\zs{\varepsilon}(d,q)=\cR^{*}_\zs{\varepsilon}(\wh{S}_\zs{d},S)\,,
\end{equation}
where the risk $\cR^{*}_\zs{\varepsilon}(\wh{S},S)$ is defined in \eqref{sec:robust-risks} and $d$ is an integer number (maybe random)
 from the set $\{1,\ldots,\iota\}$.
In this case the cost function \eqref{sec:Mo.5}
has the following form.
\begin{equation}\label{sec:NumSgn.4}
J_\varepsilon(d)=\sum_\zs{j=1}^{d} \, \wh{\theta}^2_\zs{j,\varepsilon} -2 \sum_\zs{j=1}^{d}\,\wt{\theta}_\zs{j,\varepsilon} + 
\delta\,\wh{P}_\zs{\varepsilon}(\lambda)
\,.
\end{equation}
So,  for this problem  the LSE model selection procedure
is defined as
\begin{equation}\label{sec:NumSgn.5-+-0}
\wh{q}_\zs{\varepsilon}
= \mbox{argmin}_\zs{1\le d\le \iota} J_\zs{\varepsilon}(d)\,.
\end{equation}

Note that Theorem \ref{Th.sec:Mrs.2} implies that
the robust risks of  the procedure \eqref{sec:Mo.9} with $\vert\Lambda\vert_\zs{*}\le 1/\varepsilon$, for any  $0 <\delta< 1/6$, satisfy the following oracle inequality
\begin{equation}\label{sec:NumSgn.5}
\D_\zs{\varepsilon}(\wh{q}_\zs{\varepsilon}\,,\,q)
\leq\frac{1+3\delta}{1-3\delta} \min_\zs{1\le d\le \iota}
 \D_\zs{\varepsilon}(d,q)
 +\ve^{2}\,
\frac{\U^{*}_\zs{\ve}(S)}{\delta}
\,,
\end{equation}
where the last term satisfies the property  \eqref{sec:Mrs.7-25.3}.

\bigskip

\bigskip

\section{Simulations}\label{sec:Siml}

In this section we report the results of a Monte Carlo experiment to assess the performance of the proposed model selection procedure \eqref{sec:Mo.9}. 
In \eqref{sec:In.1} we chose 
\begin{equation}\label{sec:Siml.0}
S(t)=\sum^{10}_\zs{j=1}\,\frac{j}{j+1}\,\phi_\zs{j}(t)
\,,
\end{equation}
with $\phi_\zs{j}(t)=\sqrt{2}\sin(2\pi l_\zs{j}t)$, $l_\zs{j}=[\sqrt{j}]j$. We simulate the model
$$
\d y_\zs{t} = S(t) \d t + \varepsilon\d w_\zs{t}\,.
$$

The frequency of observations  per period equals $p=100000$.
We use the weight sequence 
as proposed in \cite{GaltchoukPergamenshchikov2009a}
for a discrete time model:
$k^{*}=100+\sqrt{\vert\ln \varepsilon\vert}$ and  $m=[\vert \ln\varepsilon\vert^{2}]$. We calculated the empirical quadratic risk defined as 

$$
\overline{\R}=
\frac{1}{p}
\sum^{p}_\zs{j=1}\,
\wh{\E}\,\left(\wt{S}_\zs{\varepsilon}(u_\zs{j})-S(u_\zs{j})\right)^{2}\,,\quad u_\zs{j}=j/p\,,
$$
and the relative quadratic risk
$$
\overline{\R}_\zs{*}=\overline{\R}/\|S\|^{2}_\zs{p}
\quad\mbox{and}\quad
\|S\|^{2}_\zs{p}=\frac{1}{p}\,
\sum^{p}_\zs{j=1}\,S^{2}(u_\zs{j})
\,.
$$
\noindent 
The expectations  was taken as an average over 
 $N=10000$ replications, i.e.
$$
\wh{\E}\,\left(\wt{S}_\zs{\varepsilon}(\cdot)-S(\cdot)\right)^{2}=
\frac{1}{N}\,
\sum^{N}_\zs{l=1}\,
\left(\wt{S}^{l}_\zs{\varepsilon}(\cdot)-S(\cdot)\right)^{2}\,.
$$
 We used 
the cost function 
 with $\delta=(3+\vert\ln\varepsilon\vert)^{-2}$.

\bigskip

\bigskip

\bigskip

\centerline{ \underline{Table} : Empirical risks }

\begin{center}
    \begin{tabular}{|l|c|r|}
                                  \hline
                                 $ \varepsilon$ & $\mathbf{\overline{R}}$  & $\mathbf{\overline{R_*}}$  \\
                                  \hline
                                  $1/\sqrt{20}$ &0.0158  & 0.307  \\
                                  
                                  \hline
                                  $1/\sqrt{100}$ &0.0113  & 0.059  \\
                                  \hline
                                  $1/\sqrt{200}$ &0.0076  & 0.04  \\
                                  \hline
                                  $1/\sqrt{1000}$ & 0.0035 &0.0185  \\
                                  \hline

\end{tabular}

\end{center}

\newpage

In the following graphics the dashed line is  the model selection procedure \eqref{sec:Mo.9},
the continuous line is the function \eqref{sec:Siml.0} and the bold line is the corresponding observations \eqref{sec:In.1}.

\begin{center}

\includegraphics[scale=0.4]{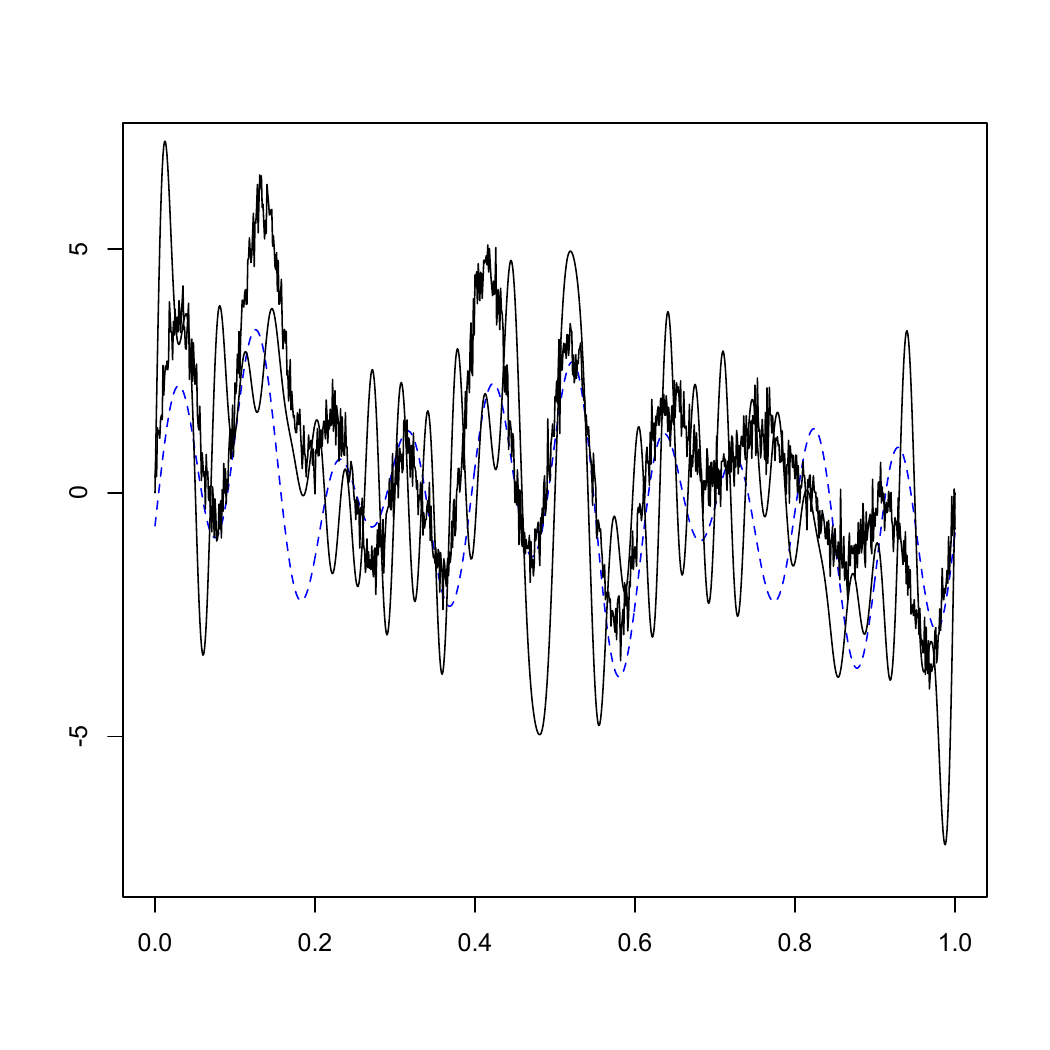}

\textbf{Figure 1: $\varepsilon=1/\sqrt{20}$}

\end{center}

\bigskip

\bigskip

\begin{center}

\includegraphics[scale=0.4]{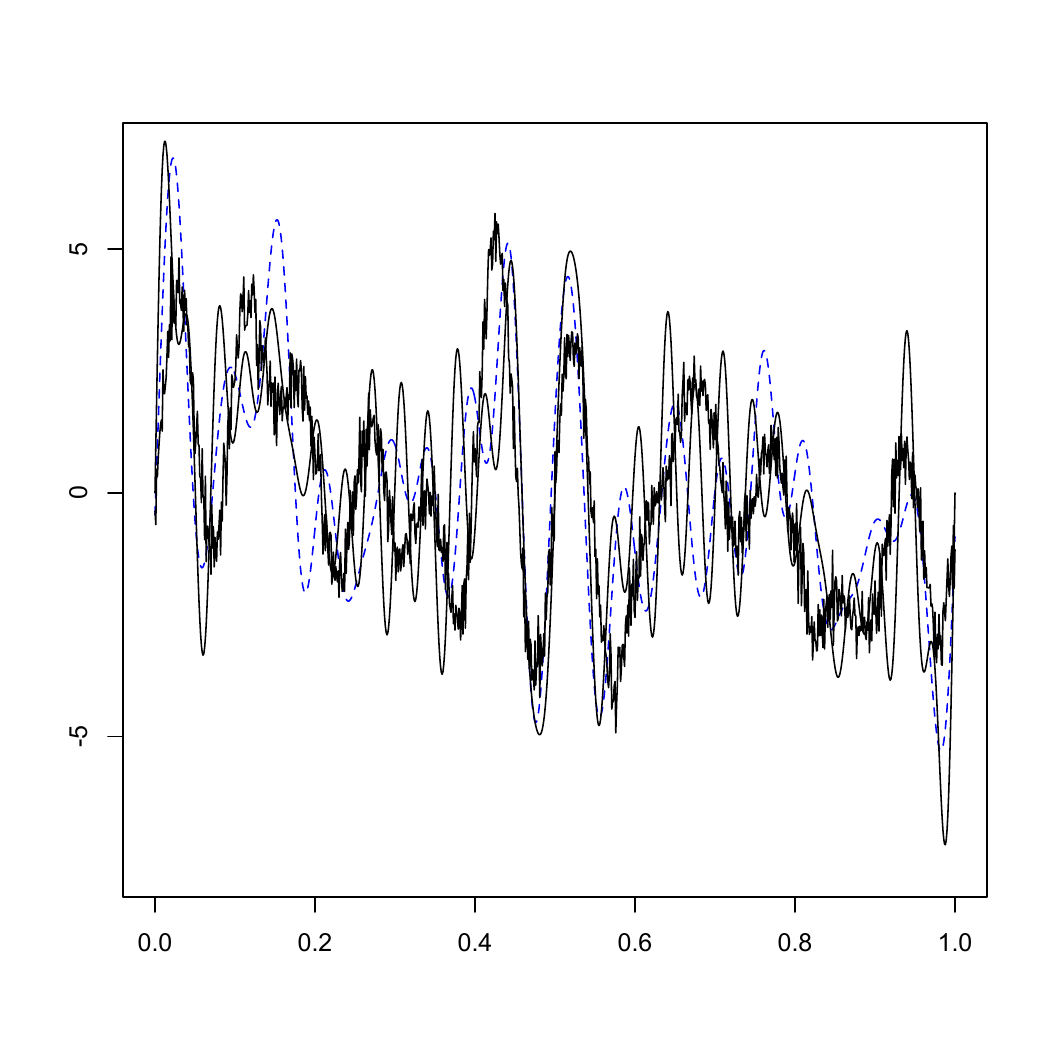}

\textbf{Figure 2: $\varepsilon=1/\sqrt{100}$}
\end{center}

\bigskip

\begin{center}

\includegraphics[scale=0.4]{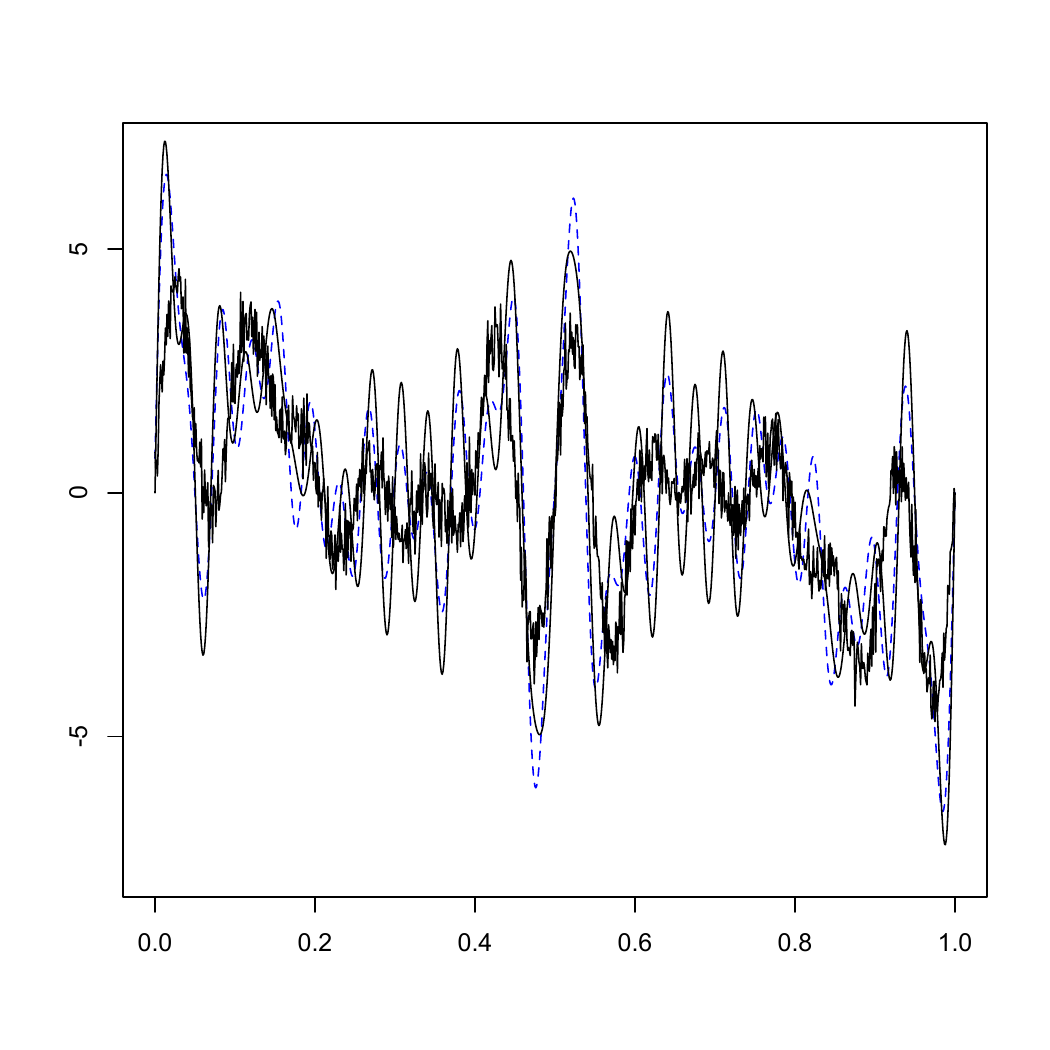}

\textbf{Figure 3: $\varepsilon=1/\sqrt{200}$}
\end{center}

\bigskip

\bigskip

\begin{center}

\includegraphics[scale=0.4]{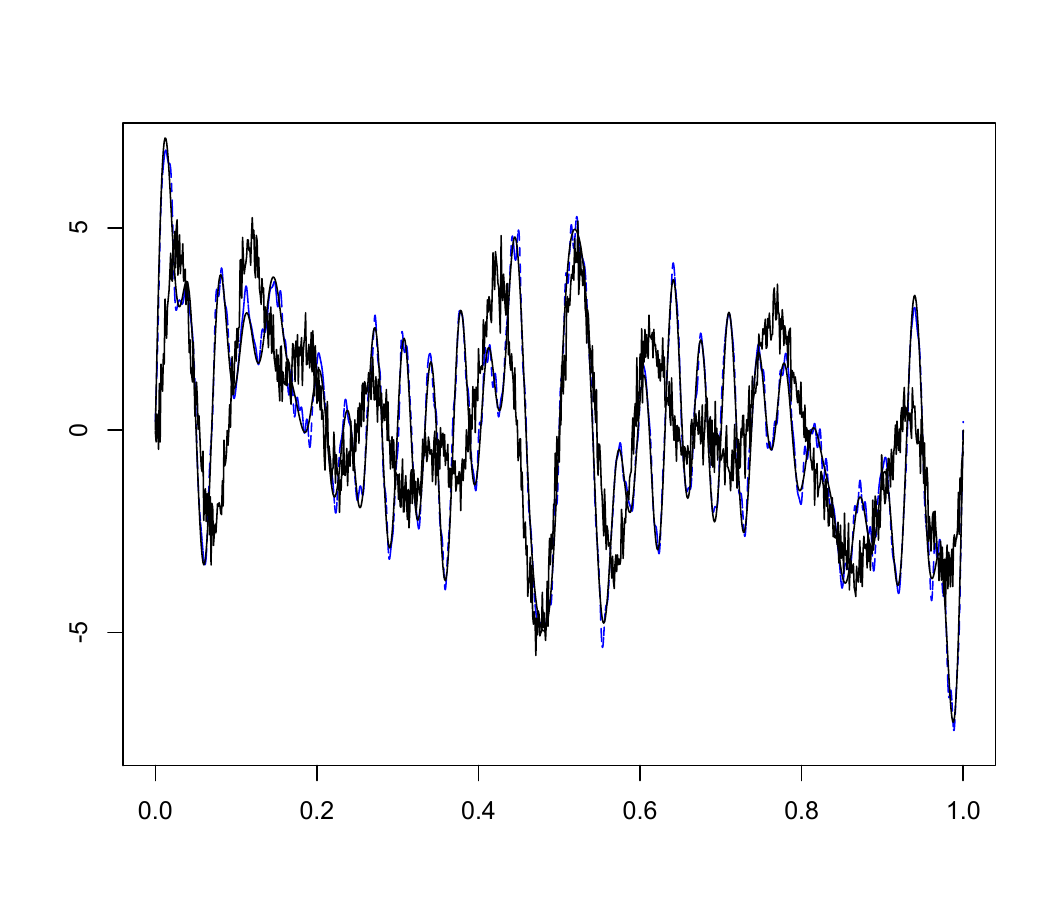}

\textbf{Figure 4: $\varepsilon=1/\sqrt{1000}$}
\end{center}

\bigskip

\bigskip

To estimate the signals number  $q$ we use two procedures. The first
$\wh{q}_\zs{1}$ is  \eqref{sec:NumSgn.5} with $\iota=[\ln \varepsilon^{-2}]$.
The second $\wh{q}_\zs{2}$ is defined through the shrinkage approach
for
the model selection procedure \eqref{sec:Siml.0}.
$$
\wh{q}_\zs{2}=\inf\{j\ge 1: \vert \wh{\theta}_\zs{j}\vert \le \c^{*}_\zs{\varepsilon}\}\,,
\quad \c^{*}_\zs{\varepsilon}=\, \varepsilon\sqrt{\vert \log\varepsilon \vert}\,.
$$

\bigskip

\bigskip

\centerline{\bf  \underline{Table} : Estimation of the number signals}

\begin{center}
    \begin{tabular}{|l|c|r|}
                                  \hline
                                  $\varepsilon$ & $\hat{q}_\zs{1}$& $\hat{q}_\zs{2}$  \\
                                  \hline
                                 $1/\sqrt{20}$  & 6 & 5  \\
                                  \hline
                                 $1/\sqrt{100}$  & 8 & 7 \\
                                  \hline
                                 $1/\sqrt{200}$ & 9 & 7   \\
                                  \hline
                                 $1/\sqrt{1000}$ & 10  & 9  \\
                                  \hline
                                 
\end{tabular}

\end{center}

\bigskip

\begin{remark}\label{re.sec:Siml.1}
From the simulation we can conclude 
that the LSE procedure \eqref{sec:NumSgn.5} is more appropriate
than shrinkage method
 for
such number detection problem. 
\end{remark}

\bigskip

\bigskip

\bigskip
\section{Proofs}\label{sec:Pr}
\subsection{Proof of Proposition \ref{Pr.sec: MaPr.1}}

First note that
\begin{equation}
\label{sec: MaPr.3}
B^{2}_\zs{2,\varepsilon}(u)\le 2\,\wt{\xi}^{2}_\zs{1}
+
2B^{2}_\zs{2,\varepsilon}(u')\,,
\end{equation}
where $u'=(0,u_\zs{2},\ldots,u_\zs{n})\in\bbr^{n}$. 
It should be noted that 
$$
\E\,\wt{\xi}^{2}_\zs{1}\le \E\,\xi^{4}_\zs{1}\le
8\left(
\varrho^{4}_\zs{1}\E\,w^{4}_\zs{1}
+
\varrho^{4}_\zs{1}\E\,z^{4}_\zs{1}
\right)=
8\left(
3\varrho^{4}_\zs{1}
+
\varrho^{4}_\zs{1}\E\,z^{4}_\zs{1}
\right)
\,.
$$

To study the last term in the right hand side of the inequality \eqref{sec: MaPr.3} we set
for any function $f$ from $\L_\zs{2}[0,1]$
$$
\check{I}_\zs{t}(f)=\int^{t}_\zs{0}\,f(s)\d\check{\xi}_\zs{s}
\quad\mbox{and}\quad
\wt{I}_\zs{t}(f)=\check{I}^{2}_\zs{t}(f)-\E\,\check{I}^{2}_\zs{t}(f)
\,.
$$
Note  that for $j\ge 2$ we define the random variables $\wt{\xi}_\zs{j}=\wt{I}_\zs{1}(\phi_\zs{j})$. So,
$$
B^{2}_\zs{2,\varepsilon}(u')=\sum^{n}_\zs{j=2}\,u_\zs{j}\,\wt{I}_\zs{1}(\phi_\zs{j})
=:D_\zs{1}(u)
\,.
$$
By the It\^o formula we can write that for any function $f$ from $\L_\zs{2}[0,1]$
$$
\d \wt{I}_\zs{t}(f)=2\check{I}_\zs{t-}(f)  f(t) \d \check{I}_\zs{t}(f)
+ \varrho^{2}_\zs{2}\,f^{2}(t)\,\d \check{m}_\zs{t}\,,
$$
where $\check{m}_\zs{t}=h^{2}_\zs{\varepsilon}*(\mu-\wt{\mu})_\zs{t}$. So, taking into account that
$$
\d \check{I}_\zs{t}(f)=\varrho_\zs{1}\,f(t)\d w_\zs{t}+
\varrho_\zs{2}\,f(t)\d \check{z}_\zs{t}
\,,
$$
we obtain that
$$
\d \wt{I}_\zs{t}(f)=2\varrho_\zs{1}\,
\check{I}_\zs{t}(f)f(t)\d w_\zs{t}
+
2\varrho_\zs{2}\,
\check{I}_\zs{t-}(f)f(t)\d \check{z}_\zs{t}
+
\varrho^{2}_\zs{2}\,f^{2}(t)\,\d \check{m}_\zs{t}\,.
$$
So, setting
$$
V_\zs{t}=\sum^{n}_\zs{j=2}\, u_\zs{j}\,\check{I}_\zs{t}(\phi_\zs{j})\phi_\zs{j}(t)
\quad\mbox{and}\quad
\Psi_\zs{t}=
\sum^{n}_\zs{j=2}u_\zs{j}\,\phi^{2}_\zs{j}(t)
\,,
$$
we obtain that
$$
\d D_\zs{t}=2\varrho_\zs{1}\,V_\zs{t}\,\d w_\zs{t}
+2\varrho_\zs{2}\,V_\zs{t-}\,\d \check{z}_\zs{t}
+
\varrho^{2}_\zs{2}\,\Psi_\zs{t}\,\d \check{m}_\zs{t}
\,.
$$
So, we obtain that
\begin{equation}
\label{sec: MaPr.4}
D^{2}_\zs{1}\le 12 \varrho^{2}_\zs{1}\,
\left(\int^{1}_\zs{0}\,V_\zs{t}\,\d w_\zs{t}\right)^{2}
+
12 \varrho^{2}_\zs{2}\,M^{2}_\zs{1}
+
3
\varrho^{4}_\zs{2}\,
\left(\int^{1}_\zs{0}\,\Psi_\zs{t-}\,\d \check{m}_\zs{t}\right)^{2}
\,,
\end{equation}
where $M_\zs{t}=\int^{t}_\zs{0}\,V_\zs{s-}(u)\,\d \check{z}_\zs{s}$.
Moreover, 
taking into account that for any $f$, $g$ from $\L_\zs{2}[0,1]$
$$
\E\,\check{I}_\zs{t}(f)\,\check{I}_\zs{t}(g)
=\check{\varkappa}_\zs{Q}\,\int^{t}_\zs{0}\,f(s)g(s)\,\d s\,,
$$
we get
\begin{align*}
2\int^{1}_\zs{0}\,\E\,V^{2}_\zs{t}\,\d t&=
2\sum^{n}_\zs{i,j=2}u_\zs{i}\,u_\zs{j}\,\int^{1}_\zs{0}\,\phi_\zs{i}(t)\phi_\zs{j}(t)\,
\E\,\check{I}_\zs{t}(\phi_\zs{i})\,\check{I}_\zs{t}(\phi_\zs{j})
\,\d t\\[2mm]
&=\check{\varkappa}_\zs{Q}\,
\sum^{n}_\zs{i,j=2}u_\zs{i}\,u_\zs{j}\,\left(\int^{1}_\zs{0}\,\phi_\zs{i}(t)\phi_\zs{j}(t)\,
\,\d t\right)^{2}
=\check{\varkappa}_\zs{Q}\,.
\end{align*}
Thus,
$$
2\E\left(\int^{1}_\zs{0}\,V_\zs{t}\,\d w_\zs{t}\right)^{2}
=\check{\varkappa}_\zs{Q}\,.
$$
Now, to estimate the second term in the inequality \eqref{sec: MaPr.4} note that
in view of the inequality \eqref{Novikov++}
for any bounded function $f$ and any $0\le t\le 1$
\begin{align*}
\E\,\check{I}^{4}_\zs{t}(f)&\le 8\varrho^{4}_\zs{1}\,\E\left( \int^{t}_\zs{0}\,f(s)\d w_\zs{s}\right)^{4}
+
 8\varrho^{4}_\zs{2}\,\E\left( \int^{t}_\zs{0}\,f(s)\d \check{z}_\zs{s}\right)^{4}
\\[2mm]
&\le 24\varrho^{4}_\zs{1}\int^{1}_\zs{0}f^{2}(t)\d t
+ C^{*}_\zs{4}\left( \left(\Pi(h^{2}_\zs{\varepsilon})\,\int^{1}_\zs{0}f^{2}(t)\d t\right)^{2}
+
\Pi(h^{4}_\zs{\varepsilon})\,\int^{1}_\zs{0}f^{4}(t)\d t
 \right)\,,
\end{align*}
i.e.
$$
\sup_\zs{0\le t\le 1}\,
\E\,\check{I}^{4}_\zs{t}(f)<\infty\,.
$$
Now it is easy to see that  through the H\"older inequality  the term $V_\zs{t}$ can be estimated as
$$
\sup_\zs{0\le t\le 1}
\E\,V^{4}_\zs{t}\,<\,\infty
\,.
$$ 
From here and the inequality \eqref{Novikov++} it follows that
$$
\sup_\zs{0\le t\le 1}\,\E\,M^{4}_\zs{t}\le C^{*}_\zs{4}\,\left(\left(\Pi(h^{2}_\zs{\varepsilon})\right)^{2}+\Pi(h^{4}_\zs{\varepsilon})\right)
\,
\int^{1}_\zs{0}\,\E\,V^{4}_\zs{t}\,\d t
<\infty
$$
and, therefore,
$$
\int^{1}_\zs{0}\,\E\, M^{2}_\zs{t} V^{2}_\zs{t}\,\d t
\le
\sup_\zs{0\le t\le 1}\,\left( \E\,M^{4}_\zs{t}\right)^{1/2}\, 
\left( \int^{1}_\zs{0}\,\E\,V^{4}_\zs{t}\,\d t
\right)^{1/2}
\,<\,\infty\,.
$$
This implies that
$$
\E\,\int^{1}_\zs{0}\,M_\zs{t-}\,\d M_\zs{t}=0\,.
$$
Thus, the It\^o formula implies
$$
2\E\,M^{2}_\zs{1}=\E\sum_\zs{0\le t\le 1}\,(\Delta M_\zs{t})^{2}
=2\Pi(h^{2}_\zs{\varepsilon})\,\int^{1}_\zs{0}\,\E\,V^{2}_\zs{t}\d t
= \Pi(h^{2}_\zs{\varepsilon})\,
\check{\varkappa}_\zs{Q}\,.
$$
In the same way we calculate
\begin{align*}
\varrho^{2}_\zs{2}\,\E\,
\left(\int^{1}_\zs{0}\,\Psi_\zs{t-}\d \check{m}_\zs{t}\right)^{2}
&=\varrho^{2}_\zs{2}\E\,\sum_\zs{0\le t\le 1}\,\left( \Delta \check{m}_\zs{t}\right)^{2}\,\Psi^{2}_\zs{t-}\\[2mm]
&=\varrho^{2}_\zs{2}\Pi(h^{4}_\zs{\varepsilon})\,\int^{1}_\zs{0}\,\Psi^{2}_\zs{t}\d t
\le \Pi(x^{2})(a/\varepsilon)^{2}\,\left(\phi^{*}\right)^{4}\,\#(u)\,.
\end{align*}
So, taking into account that $\Pi(x^{2})=1$,
 we obtain that
$$
\E\,D^{2}_\zs{1}\le 6\varrho^{2}_\zs{1}\,\check{\varkappa}_\zs{Q}
+3\varrho^{2}_\zs{2}\,
\left(2\check{\varkappa}_\zs{Q} +(\phi^{*})^{4}\right)
\le 6\,\varkappa^{2}_\zs{Q}
+
3\varrho^{2}_\zs{2}
\,
(\phi^{*})^{4}
\,.
$$
Similarly we obtain that
$$
\E\wt{\xi}^{2}_\zs{1}\le 
 6
 \varkappa^{2}_\zs{Q}
+
3\varrho^{4}_\zs{2}\,\Pi(x^{4})
\,.
$$
This implies the upper bound  \eqref{sec: MaPr.2}. \endproof

\bigskip

\bigskip

\subsection{Proof of Thoerem \ref{Th.sec:Mrs.1} }

First note, that
 we can rewrite the empirical squared error in \eqref{sec:Mo.3} as follows
\begin{equation}\label{sec:Pr.1}
\Er_\varepsilon(\lambda) = J_\varepsilon(\lambda) + 2 \sum_\zs{j=1}^{n} \lambda(j) \check{\theta}_\zs{j,\varepsilon}+ 
\Vert S\Vert^2-\delta 
\wh{P}_\varepsilon(\lambda),
\end{equation}
where $\check{\theta}_\zs{j,\varepsilon}=\wt{\theta}_\zs{j,\varepsilon}-\theta_\zs{j}\wh{\theta}_\zs{j,\varepsilon}$. 
Now using the definition of $\wt{\theta}_\zs{j,\varepsilon}$ in \eqref{sec:Mo.4}
we obtain that
$$
\check{\theta}_\zs{j,\varepsilon}=\varepsilon\theta_\zs{j}\overline{\xi}_\zs{j} +\varepsilon^2 \wt{\xi}_\zs{j}  +\varepsilon^2 
(\E_\zs{Q}\,\overline{\xi}^{2}_\zs{j,\varepsilon}-\check{\varkappa}_\zs{Q})
 +
\varepsilon^2 (\check{\varkappa}_\zs{Q} -  \wh{\varkappa}_\zs{\varepsilon}) 
\,.
$$
Where $\wt{\xi}_\zs{j}= \overline{\xi}^{2}_\zs{j}- \E_\zs{Q}\,\overline{\xi}^{2}_\zs{j}$ and $\overline{\xi}_\zs{j}=  \check{I}_\zs{1}(\phi_\zs{j})$.Setting
\begin{equation}\label{sec:Pr.2}
M_\zs{\varepsilon}(\lambda) = \varepsilon\sum_\zs{j=1}^{n} \lambda(j)\theta_\zs{j} \overline{\xi}_\zs{j}
\quad\mbox{and}\quad
L(\lambda)=
\sum_\zs{j=1}^{n} \lambda(j)
\,,
\end{equation}
we can rewrite \eqref{sec:Pr.1} as
\begin{align}\nonumber
\Er_\varepsilon(\lambda) & =  J_\varepsilon(\lambda) + 2\varepsilon^2 (\check{\varkappa}_\zs{Q}-  \wh{\varkappa}_\zs{\varepsilon} )\,L(\lambda)+ 2 
M_\zs{\varepsilon}(\lambda)+2\varepsilon^2  B_\zs{1,\varepsilon}(\lambda)\\  \label{sec:Pr.3}
& +  2 \varepsilon\,\sqrt{P_\zs{\varepsilon}(\lambda)} \frac{B_\zs{2,\varepsilon}(u_\zs{\lambda})}{\sqrt{\check{\varkappa}_\zs{Q} }} + \Vert S\Vert^2-\delta 
\wh{P}_\varepsilon(\lambda),
\end{align}
where $u_\zs{\lambda}=\lambda/|\lambda|$, the exact penalization is defined in
\eqref{sec:Mo.6+1}
 and the functions $B_\zs{1,\varepsilon}(\cdot)$ and $B_\zs{2,\varepsilon}(\cdot)$ are defined in \eqref{sec: MaPr.1}.
 It should be noted that for the truncation parameter \eqref{sec:Mo.2++}
 the bound
 \eqref{sec: MaPr.2}
implies 
 \begin{equation}
\label{sec: Pr.1+}
\sup_\zs{\lambda\in\Lambda}\,
\E_\zs{Q}\,
\left\vert B^{2}_\zs{2,\varepsilon}(u_\zs{\lambda}) \right\vert
\le 
U_\zs{Q}
+
6\check{\varkappa}_\zs{Q}\,
\left(
\frac{\overline{a}}{\varepsilon}
\right)^{2}
\,\vert\Lambda\vert_\zs{*}\,(\phi^{*})^{4}
=
U_\zs{1,Q}
\,,
\end{equation}
where $
U_\zs{1,Q}
=
U_\zs{Q}
+
6\check{\varkappa}_\zs{Q}\,(\phi^{*})^{4}$.

Let $\lambda_0= (\lambda_0(j))_\zs{1\le j\le\,n}$ be a fixed sequence in $\Lambda$ and $\wh{\lambda}$ be as in \eqref{sec:Mo.9}.
Substituting $\lambda_0$ and $\wh{\lambda}$ in the equation \eqref{sec:Pr.3}, we obtain
\begin{align}
\Er_\varepsilon(\wh{\lambda})-\Er_\varepsilon(\lambda_0)& = J(\wh{\lambda})-J(\lambda_0)+
2\varepsilon^2 (\check{\varkappa}_\zs{Q}-  \wh{\varkappa}_\zs{\varepsilon} )\,L(\varpi)\nonumber\\[2mm]
&+ 2\varepsilon^2 B_\zs{1,\varepsilon}(\varpi)+2 M_\zs{\varepsilon}(\varpi)\nonumber\\[2mm]
& + 2\varepsilon \sqrt{P_\zs{\varepsilon}(\wh{\lambda})} \frac{B_\zs{2,\varepsilon}(\wh{u})}{\sqrt{\check{\varkappa}_\zs{Q}}}- 2\varepsilon \sqrt{P_\zs{\varepsilon}(\lambda_0)}
 \frac{B_\zs{2,\varepsilon}(u_\zs{0})}{\sqrt{\check{\varkappa}_\zs{Q}}}\nonumber \\[2mm]\label{sec:Pr.4}
& -  \delta  \wh{P}_\zs{\varepsilon}(\wh{\lambda})+
\delta \wh{P}_\zs{\varepsilon}(\lambda_0)
\end{align}
where $\varpi= \wh{\lambda} - \lambda_\zs{0}$, $\wh{u} = u_\zs{\wh{\lambda}}$ and $u_\zs{0} = u_\zs{\lambda_\zs{0}}$. Note that by \eqref{sec:Mo.2}
$$ 
|L(\varpi)| \le\,L(\hat \lambda) + L(\lambda) \leq 2\vert\Lambda\vert_\zs{*}
\,. 
$$
The inequality
\begin{equation}\label{sec:Pr.5}
2|ab| \leq \delta a^2 + \delta^{-1} b^2
\end{equation}
implies that for any $\lambda\in\Lambda$
$$
2 \varepsilon \sqrt{P_\zs{\varepsilon}(\lambda)} \frac{|B_\zs{2,\varepsilon}(u_\zs{\lambda})|}{\sqrt{\check{\varkappa}_\zs{Q}}} \le\, \delta P_\zs{\varepsilon}(\lambda) + \varepsilon^2
 \frac{B^2_\zs{2,\varepsilon}(u_\zs{\lambda})}{\delta\check{\varkappa}_\zs{Q}}.
$$
From the bound \eqref{sec: MaPr.1+1} it follows that
for $0 < \delta < 1$
\begin{align*}
\Er_\varepsilon(\hat \lambda) & \le\,\Er_\varepsilon(\lambda_0) +2 M_\zs{\varepsilon}(\varpi)+ 
2\varepsilon^2\frac{B^*_\zs{2,\varepsilon}}{\delta\check{\varkappa}_\zs{Q}}
+2\varepsilon^{2}\,\check{\varkappa}_\zs{Q} \\[2mm]
& +  \varepsilon^2 |\wh{\varkappa} -\check{\varkappa}_\zs{Q}| ( |\wh{\lambda}|^2 + |\lambda_0|^2+4\vert\Lambda\vert_\zs{*})+ 2  \delta P_\varepsilon(\lambda_0)\,,
\end{align*}
where $B^*_\zs{2,\varepsilon} = \sup_\zs{\lambda\in\Lambda} B^2_\zs{2,\varepsilon}(u_\zs{\lambda})$. 
It should be noted that  
through
\eqref{sec: Pr.1+}
 we can estimate  this term  as
\begin{equation}\label{sec:Pr.5++}
\E_\zs{Q}\, B^*_\zs{2,\varepsilon} \leq \sum_\zs{\lambda\in\Lambda}\E_\zs{Q}\, 
B^2_\zs{2,\varepsilon} (u_\zs{\lambda}) \leq \iota 
U_\zs{1,Q}\,.
\end{equation}
Taking into account that $\sup_\zs{\lambda\in\Lambda} |\lambda|^2 \leq \vert\Lambda\vert_\zs{*}$,
we can rewrite the previous  bound as
\begin{align}\label{sec:Pr.6}
\Er_\varepsilon(\wh{\lambda}) & \le  \Er_\varepsilon(\lambda_0) +2 M_\zs{\varepsilon}(\varpi)
+2\varepsilon^2 \frac{B^*_\zs{2,\varepsilon}}{\delta\check{\varkappa}_\zs{Q}}
+2\varepsilon^{2}\,\check{\varkappa}_\zs{Q}  \nonumber\\[2mm]
& +  \frac{6\varepsilon^2 \vert\Lambda\vert_\zs{*}}{n} |\wh{\varkappa} -\check{\varkappa}_\zs{Q}| + 2  \delta P_\varepsilon(\lambda_0).
\end{align}
To estimate the second term in the right hand side of this inequality we introduce
$$
S_\zs{\upsilon} = \sum_\zs{j=1}^{n} \upsilon(j) \theta_\zs{j} \phi_\zs{j}\,,
\quad \upsilon=(\upsilon(j))_\zs{1\le j\le n}\in\bbr^{n}\,.
$$
Moreover, note  that
$$
M^2_\zs{\varepsilon} (\upsilon)
\le 2\varepsilon^{2}
\left(\upsilon^{2}(1)\,\xi^{2}_\zs{1}
+\check{I}_\zs{1}(\Phi)
\right)\,,
$$
where $\Phi(t)=\sum^{n}_\zs{j=2}\,\upsilon(j)\theta_\zs{j}\phi_\zs{j}(t)$.
Therefore, thanks to 
\eqref{sec:In.3} we obtain that
for any nonrandom $\upsilon\in\bbr^{n}$ 
\begin{equation}\label{sec:Pr.8--0}
\E M^2_\zs{\varepsilon} (\upsilon) \leq 2\check{\varkappa}_\zs{Q} \varepsilon^2 \sum_\zs{j=1}^{n} \upsilon^{2}(j) \theta^2_\zs{j} =2 \check{\varkappa}_\zs{Q}
 \varepsilon^2 ||S_\zs{\upsilon}||^2
 \,.
\end{equation}
To estimate this function for a random vector we set
$$ 
M^*_\zs{\varepsilon} = \sup_\zs{\upsilon \in \Lambda_1} \frac{M^2 (\upsilon)}{\varepsilon^2 ||S_\zs{\upsilon}||^2}
\quad\mbox{and}\quad
\Lambda_\zs{1} = \Lambda - \lambda_\zs{0}\,.
$$
So, through the inequality  \eqref{sec:Pr.5}
\begin{equation}\label{sec:Pr.10}
2 |M_\zs{\varepsilon}(\upsilon)|\leq \delta ||S_\zs{\upsilon}||^2 +\varepsilon^2 \frac{M^*_\zs{\varepsilon}}{\delta}\,.
\end{equation}
It is clear that the last term  here can be estimated as
\begin{equation}\label{sec:Pr.10++}
\E\, M^*_\zs{\varepsilon} \leq \sum_\zs{\upsilon \in \Lambda_\zs{1}}\,
 \frac{\E\, M^2_\zs{\varepsilon}(\upsilon)}{\varepsilon^{2}||S_\zs{\upsilon}||^2} \leq 2\sum_\zs{\upsilon \in \Lambda_1} \check{\varkappa}_\zs{Q}= 
 2\check{\varkappa}_\zs{Q}\,\iota\,,
\end{equation}
where $\iota = \#(\Lambda)$. 
Moreover, note that, for any $\upsilon\in\Lambda_1$,
$$
||S_\zs{\upsilon}||^2-||\wh{S}_\zs{\upsilon}||^2 = \sum_\zs{j=1}^{n} \upsilon^{2}(j) (\theta^2_\zs{j}-\wh{\theta}^2_\zs{j}) \le 2 
\vert M_\zs{\varepsilon}(\upsilon^2)\vert\,,
$$
where $\upsilon^{2} =(\upsilon^{2}(j))_\zs{1\le j\le n}$.
Taking into account now, that for any $x \in \Lambda_1$ the components $|\upsilon(j)|\leq 1$ ,  we can estimate the last term as
in \eqref{sec:Pr.8--0}, i.e.
$$
\E\, M^2_\zs{\varepsilon}(\upsilon^{2}) \leq 2\varepsilon^2 \check{\varkappa}_\zs{Q}\,
||S_\zs{\upsilon}||^2\,.
$$
Similarly,  setting
$$ 
M^*_\zs{1,\varepsilon} = \sup_\zs{\upsilon \varepsilon \Lambda_1}\,
\frac{ M^2_\zs{\varepsilon}(\upsilon^{2})}{\varepsilon^{2}||S_\zs{\upsilon}||^2}
$$
we obtain
\begin{equation}\label{sec:Pr.10+++}
\E_\zs{Q}\, M^*_\zs{1,\varepsilon} \leq 2\check{\varkappa}_\zs{Q}\,\iota\,.
\end{equation}
By the same way we find that
$$
2 |M_\zs{\varepsilon}(\upsilon^{2})|\leq \delta ||S_\zs{\upsilon}||^2 + \frac{M^*_\zs{1,\varepsilon}}{n\delta}
$$
and, for any $0<\delta<1$,
$$
||S_\zs{\upsilon}||^2 \leq \frac{||\wh{S}_\zs{\upsilon}||^2}{1-\delta} +  \frac{\varepsilon^{2} M^*_\zs{1,\varepsilon}}{ \delta (1-\delta)}
\,.
$$
So, from \eqref{sec:Pr.10} 
we get
$$
2 M(\upsilon) \leq \frac{\delta ||\wh{S}_\zs{\upsilon}||^2}{1-\delta} +  \frac{\varepsilon^{2}(M^*_\zs{\varepsilon}+M^*_\zs{1,\varepsilon})}{ \delta (1-\delta)}
\,.
$$
Therefore, taking into account that $\Vert\wh{S}_\zs{\varpi}\Vert^{2}\le 2\,(\Er_\varepsilon(\wh{\lambda})+\Er_\varepsilon(\lambda_0))$, 
the term $M_\zs{\varepsilon}(\varpi)$ can be estimated as
$$
2 M_\zs{\varepsilon}(\varpi) \leq \frac{2\delta(\Er_\varepsilon(\wh{\lambda})+\Er_\varepsilon(\lambda_0))}{1-\delta} + 
\, \frac{\varepsilon^{2}(M^*_\zs{\varepsilon}+M^*_\zs{1,\varepsilon})}{ \delta (1-\delta)}.
$$
Using this bound in  \eqref{sec:Pr.6} 
we obtain 
\begin{align*}
\Er_n(\wh{\lambda}) & \le \frac{1+\delta}{1-3\delta} \Er_\varepsilon(\lambda_0) 
+ \frac{\varepsilon^2(M^*_\zs{\varepsilon}+M^*_\zs{1,\varepsilon})}{ \delta (1-3\delta)}
+\frac{2\varepsilon^2 B^*_\zs{2,\varepsilon}}{\delta(1-3\delta)\check{\varkappa}_\zs{Q}} \\[2mm]
& 
+\frac{2\varepsilon^{2}\,\check{\varkappa}_\zs{Q}}{1-3\delta}
+  \frac{6\varepsilon^{2}\,\vert\Lambda\vert_\zs{*}}{(1-3\delta)} |\wh{\varkappa} -\check{\varkappa}_\zs{Q}| + \frac{2\delta}{(1-3\delta)} P_\zs{\varepsilon}(\lambda_0).
\end{align*}
Moreover, for $0<\delta<1/6$ we can rewrite this inequality as
\begin{align*}
\Er_n(\wh{\lambda}) & \le \frac{1+\delta}{1-3\delta} \Er_\varepsilon(\lambda_0) 
+ \frac{2\varepsilon^2(M^*_\zs{\varepsilon}+M^*_\zs{1,\varepsilon})}{ \delta}
+\frac{4\varepsilon^2 B^*_\zs{2,\varepsilon}}{\delta\check{\varkappa}_\zs{Q}} \\[2mm]
& 
+4\varepsilon^{2}\,\check{\varkappa}_\zs{Q}
+  12\varepsilon^{2}\,\vert\Lambda\vert_\zs{*} |\wh{\varkappa} -\check{\varkappa}_\zs{Q}| + 
4\delta\,
P_\zs{\varepsilon}(\lambda_0)\,.
\end{align*}
Using here the bounds \eqref{sec:Pr.5++}, \eqref{sec:Pr.10++}, \eqref{sec:Pr.10+++}
and taking into account that $\check{\varkappa}_\zs{Q}\le \varkappa_\zs{Q}$
we obtain that
\begin{align*}
 \cR(\wh{S}_*,S)
& \le \frac{1+\delta}{1-3\delta} \cR(\wh{S}_\zs{\lambda_0},S)
+ \frac{4\varepsilon^{2}\,\varkappa_\zs{Q}(2\iota+\delta)}{\delta}
+\frac{4\varepsilon^2 
U_\zs{1,Q}\iota}{\delta\check{\varkappa}_\zs{Q}} \\[2mm]
& 
+  12\varepsilon^{2}\,\vert\Lambda\vert_\zs{*} \E_\zs{Q}\,|\wh{\varkappa} -\check{\varkappa}_\zs{Q}| + 
\frac{2\delta}{1-3\delta}\,
P_\zs{\varepsilon}(\lambda_0)\,.
\end{align*}
Now, from  Lemma~\ref{Le.sec:A.1-06-11-01} 
it follows that
\begin{align*}
 \cR(\wh{S}_*,S)
& \le \frac{1+3\delta}{1-3\delta} \cR(\wh{S}_\zs{\lambda_0},S)
+ \frac{4\varepsilon^{2}\varkappa_\zs{Q}(2\iota+\delta)}{\delta}
+\frac{4\varepsilon^2 
U_\zs{1,Q}\iota}{\delta\check{\varkappa}_\zs{Q}} \\[2mm]
& 
+  12\varepsilon^{2}\,\vert\Lambda\vert_\zs{*} \E_\zs{Q}\,|\wh{\varkappa} -\check{\varkappa}_\zs{Q}| + 
\varepsilon^{2}\frac{2\delta}{1-3\delta}\,\varkappa_\zs{Q}\,.
\end{align*}
Taking into account here 
 that $2\delta/(1-3\delta)\le 1$
for $0<\delta<1/6$ 
and using the function
 \eqref{sec:OIn.1}
we obtain
 the inequality
\eqref{sec:OI.1} for some constant $\l_\zs{*}>0$ which depends on $\Pi(x^{4})$. 
Hence Theorem \ref{Th.sec:Mrs.1}. 
\fdem

\subsection{Proof of Proposition \ref{Pr.sec: MaPr.3} }

We use here the same method as in \cite{KonevPergamenshchikov2009a}.
First, note that from the definitions \eqref{sec:In.8} and \eqref{sec:Mo.4-2-31-3}
we obtain
\begin{equation}\label{sec:Mo.1-1-04}
\wh{\tau}_\zs{j,\varepsilon}= \tau_\zs{j}+
\varepsilon\,
\eta_\zs{j}\,,
\end{equation}
where  
$$
\tau_\zs{j}=
\int^{1}_\zs{0}\,S(t)\,\Tr_\zs{j}(t)\d t
\quad\mbox{and}\quad
\eta_\zs{j}=
\int^{1}_\zs{0}\,\Tr_\zs{j}(t)\,\d \check{\xi}_\zs{t}
\,.
$$
So, we have
\begin{equation}\label{sec:Mo.1-0-1-04}
\wh{\varkappa}_\zs{\varepsilon}=
\sum^n_\zs{j=[1/\varepsilon]+1}\,\tau^2_\zs{j}
+
2
\check{M}_\zs{\varepsilon}
+
\varepsilon^{2}
\sum^n_\zs{j=[1/\varepsilon]+1}\,\eta^{2}_\zs{j}
\,,
\end{equation}
where $\check{M}_\zs{\varepsilon}=\varepsilon\,\sum^n_\zs{j=[1/\varepsilon]+1}\,\tau_\zs{j}\,\eta_\zs{j}$.
Note that for the continuously 
differentiable functions (see, for example, Lemma A.6 in \cite{KonevPergamenshchikov2009a})
the Fourrier coefficients $(\tau_\zs{j})$
for any $n\ge 1$ satisfy the following inequality
\begin{equation}\label{sec:Mo.2-1-04}
\sum^{\infty}_\zs{j=[1/\varepsilon]+1}\,\tau^2_\zs{j}
\le 
4\varepsilon\left(\int^{1}_\zs{0}\vert\dot{S}(t)\vert \d t\right)^{2}
\le 4\varepsilon
\Vert\dot{S}\Vert^{2}
\,.
\end{equation}
We recall, that $\dot{S}$ is the derivative of $S$. 
The term $\check{M}_\zs{\varepsilon}$
can be estimated by the same way as in
\eqref{sec:Pr.8--0}, i.e.
$$
\E_\zs{Q}\,\check{M}^{2}_\zs{\varepsilon}\le \check{\varkappa}_\zs{Q}\varepsilon^{2}\,
\sum^{n}_\zs{j=[1/\varepsilon]+1}\,\tau^{2}_\zs{j}
\le 
4\varepsilon^{3}\check{\varkappa}_\zs{Q}\Vert\dot{S}\Vert^{2}\,.
$$
Moreover, taking into account that
for $j\ge 2$ the expectation $\E\,\eta^{2}_\zs{j}=\check{\varkappa}_\zs{Q}$
we can represent the last term in \eqref{sec:Mo.1-0-1-04} as
$$
\varepsilon^{2}
\sum^n_\zs{j=[1/\varepsilon]+1}\,\eta^{2}_\zs{j}
=
\check{\varkappa}_\zs{Q}(\varepsilon^{2} n- \varepsilon^{2}[1/\varepsilon])
+ \varepsilon
\,
B_\zs{2,\varepsilon}(x')
\,,
$$
where the function $B_\zs{2,\varepsilon}(x')$ is defined in  
\eqref{sec: MaPr.1}
and $x'_\zs{j}=\varepsilon\Chi_\zs{\{1/\varepsilon<j\le 1/\varepsilon^{2}\}}$. We remind that $n=[1/\varepsilon^{2}]$.
Therefore, in view of Proposition \eqref{Pr.sec: MaPr.1}
we obtain
$$
\E_\zs{Q}
\left\vert
\varepsilon^{2}
\sum^n_\zs{j=[\sqrt{1/\varepsilon}]+1}\,\eta^2_\zs{j}
-
\check{\varkappa}_\zs{Q}
\right\vert
\le 
2\varepsilon\,
\check{\varkappa}_\zs{Q}
+
\varepsilon \sqrt{U_\zs{Q}}+\frac{\sqrt{6\check{\varkappa}_\zs{Q}}}{\vert\Lambda\vert_\zs{*}}
\,.
$$
So, we obtain the bound \eqref{sec: MaPr.5--1}. Hence 
Proposition \ref{Pr.sec: MaPr.3}\,.
\endproof

\bigskip

\bigskip

\bigskip

{\bf Acknowledgements.}  

The construction of the model selection procedures, the construction of the signals detection procedures  and numerical simulations were 
 done under financial support of the RSF Grant Number 14-49-00079 (National Research University “MPEI” 14 Krasnokazarmennaya, 111250 Moscow, Russia).
The oracles inequalities were obtained under financial support of 
 the RSF grant 17-11-01049 (National Research Tomsk State University). The analyse efficiency:  the upper and lower bounds for the robust risk, 
 the  van Trees inequalities for the impulse noises, the calculation of the Pinsker constant  were done under financial support 
     of the RFBR Grant no. 16-01-00121,
     the Ministry of Education and Science of the Russian Federation in the framework of the research project no. 2.3208.2017/4.6 
  and  the Russian Federal Professor program (project no. 1.472.2016/1.4, 
  Ministry of Education and Science).

\bigskip

\renewcommand{\theequation}{A.\arabic{equation}}
\renewcommand{\thetheorem}{A.\arabic{theorem}}
\renewcommand{\thesubsection}{A.\arabic{subsection}}
\section{Appendix}\label{sec:A}
\setcounter{equation}{0}
\setcounter{theorem}{0}

\subsection{Property of the penalty term}

\begin{lemma}\label{Le.sec:A.1-06-11-01}
Assume that Proposition \ref{Pr.sec: MaPr.0} holds. Then for any $n\ge\,1$ and $\lambda \in \Lambda$,
$$ 
P_\zs{\varepsilon}(\lambda) \leq 
\cR(\wh{S}_\zs{\lambda},S)
+
\varepsilon^{2}\,
\varkappa_\zs{Q}\,. 
$$
where the coefficient $P_\zs{\varepsilon}(\lambda)$ is defined in \eqref{sec:Mo.6+1}.
\end{lemma}
\proof In vue of the definition of $\Er_\varepsilon(\lambda)$ and the equation \eqref{sec:In.8}
  one has
\begin{align*}
\Er_\varepsilon(\lambda)&= 
\sum_\zs{j=1}^{n} \left((\lambda(j)-1) \theta_\zs{j}+\varepsilon \lambda(j)\overline{\xi}_\zs{j} \right)^{2}
+\sum^{\infty}_\zs{j=n+1}\theta^{2}_\zs{j}
\\[2mm]
&
\ge 
2\varepsilon
\sum_\zs{j=1}^{n}\,(\lambda(j)-1) \theta_\zs{j} \lambda(j)\overline{\xi}_\zs{j} 
+\varepsilon^2
\sum_\zs{j=1}^{n} \lambda^{2}(j)\overline{\xi}^{2}_\zs{j}
\,. 
\end{align*}                                     
Moreover,
using here the definition
\eqref{sec: MaPr.1} 
 we obtain 
$$
\E_\zs{Q}\, \Er_\varepsilon(\lambda) \ge\, \varepsilon^2 \sum_\zs{j=1}^{n} \lambda^2(j) 
 \E_\zs{Q}\,\overline{\xi}^{2}_\zs{j}=P_\zs{\varepsilon}(\lambda)
 - \varepsilon^2 B_\zs{1,\varepsilon}(\lambda^{2}) 
\,,
$$
where $\lambda^{2}=(\lambda^{2}(j))_\zs{1\le j\le n}$. Now,
Proposition \ref{Pr.sec: MaPr.0} implies Lemma \ref{Le.sec:A.1-06-11-01}.
\fdem

\subsection{Proof of the limit equality  \eqref{sec:Lo.12_11.1}}

First, setting $\zeta_\zs{\varepsilon}=\sum^{d}_\zs{j=1}\,\kappa^{2}_\zs{j}\,a_\zs{j}$, 
we obtain that
$$
\left\{
S_\zs{\kappa}\notin W_\zs{k,\r}
\right\}
=
\left\{
\zeta_\zs{\varepsilon}
>\r
\right\}
\,.
$$
Moreover, note that one can check directly that
$$
\lim_\zs{\varepsilon\to 0}\,
\E\,\zeta_\zs{\varepsilon}=
\lim_\zs{\varepsilon\to 0}\,
\frac{1}{v_\zs{\varepsilon}}
\sum^{d}_\zs{j=1}\,s^{*}_\zs{j}\,a_\zs{j}=\check{\r}=
(1-\check{\gamma})\r
\,.
$$
So, for sufficiently small $\varepsilon$ we obtain that
$$
\left\{
S_\zs{\kappa}\notin W_\zs{k,r}
\right\}
\subset 
\left\{
\wt{\zeta}_\zs{\varepsilon}>
\r_\zs{1}
\right\}
\,,
$$
where $\r_\zs{1}=\r\check{\gamma}/2$,
$
\wt{\zeta}_\zs{\varepsilon}=\zeta_\zs{\varepsilon}-\E\,\zeta_\zs{\varepsilon}
=v^{-1}_\zs{\varepsilon}\,\sum^{d}_\zs{j=1}\,s^{*}_\zs{j}a_\zs{j}\wt{\eta}_\zs{j}$
and
$\wt{\eta}_\zs{j}=\eta^{2}_\zs{j}-1$
Through the correlation inequality 
 (see, Proposition A.1 in \cite{GalthoukPergamenshchikov_2013})
we can get that for any $p\ge 2$ 
$$
\E\,\wt{\zeta}^{p}_\zs{\varepsilon}\le (2p)^{p/2}
\E\vert \wt{\eta}_\zs{1}\vert^{p}
\,v^{-p}_\zs{\varepsilon}\,
\left( \sum^{d}_\zs{j=1}\,
(s^{*}_\zs{j})^{2}a^{2}_\zs{j}
\right)^{p/2}
=\mbox{O}(\,v^{-\frac{p}{4k+2}}_\zs{\varepsilon})
\,,
$$
 as $\varepsilon\to 0$. Therefore,
 for any $\iota>0$
  using the Chebychev inequality for 
 $p>(4k+2)\iota$
we obtain that 
$$
v^{\iota}_\zs{\varepsilon}\P(\wt{\zeta}_\zs{\varepsilon}>\r_\zs{1})\to 0
\quad\mbox{as}\quad
\varepsilon\to 0\,.
$$
Hence the equality \eqref{sec:Lo.12_11.1}.
\endproof

\bigskip

\bigskip

\subsection{The absolute continuity of distributions for the L\'evy processes.}\label{subsec:App.5++}

In this section we study the absolute continuity for the  the L\'evy processes defined as

\begin{equation}\label{sec:App.5++.1}
  \d y_t=S(t)\d t+\d \xi_t\,,
  \quad 0\le t\le T\,,
 \end{equation}
where $S(\cdot)$ is any arbitrary nonrandom square integrated function, i.e. from $\L_\zs{2}[0,T ]$ and 
$(\xi_\zs{t})_\zs{0\le t\le T}$ is the L\'evy process of the form \eqref{sec:In.1+1}
with nonzero constants $\varrho_\zs{1}$
and $\varrho_\zs{2}$. We  denote by $\P_\zs{y}$ and $\P_\zs{\xi}$ the distributions of the processes 
$(y_\zs{t})_\zs{0\le t\le 1}$
and
$(\xi_\zs{t})_\zs{0\le t\le 1}$
on the Skorokhod space $\D[0,T]$.  Now for any $0\le t\le T$ and $(x_\zs{t})_\zs{0\le t\le T}$
from $\D[0,T]$
 we set
\begin{equation}\label{sec:App.++.1}
\Upsilon_\zs{t}(x)=
\exp\left\{\int^{t}_\zs{0}\,\frac{S(u)}{\varrho^{2}_\zs{1}}\,\d x^{c}_\zs{u}
-\,\int^{t}_\zs{0}\,
\frac{S^{2}(u)}{2\varrho^{2}_\zs{1}}\,
\d u
\right\}
\,,
\end{equation}
where $x^{c}$ is the continuous part of the process $x$ defined in \eqref{sec:App.7}.
Now we study the measures $\P_\zs{y}$ and 
$\P_\zs{\xi}$ in $\D[0,T]$.
\begin{proposition}\label{Pr.sec:App.1++}
For any  $T>0$ the measure $\P_\zs{y}\ll \P_\zs{\xi}$ 
in $\D[0,T]$
and the Radon-Nikodym derivative
is 
$$
\frac{\d\P_\zs{y}}{\d\P_\zs{\xi}}(\xi)
=\Upsilon_\zs{T}(\xi)
\,.
$$
\end{proposition}
\noindent {\bf  Proof.}
Note that to show  this proposition it suffices to check that
for any $0=t_\zs{0}<\ldots<t_\zs{n}=T$
any $b_\zs{j}\in\bbr$ for $1\le j\le n$
$$
\E\,\exp\left\{i\sum^{n}_\zs{l=1}b_\zs{j}(y_\zs{t_\zs{j}}-y_\zs{t_\zs{j-1}})\right\}
=
\E\,\exp\left\{i\sum^{n}_\zs{l=1}b_\zs{j}(\xi_\zs{t_\zs{j}}-\xi_\zs{t_\zs{j-1}})\right\}\Upsilon_\zs{T}(\xi)
\,.
$$
taking into account that the processes $(y_\zs{t})_\zs{0\le t\le T}$ and 
$(\xi_\zs{t})_\zs{0\le t\le T}$
have the independent homogeneous increments, to this end one needs to check only 
that for any $b\in\bbr$ and $0\le s<t\le T$
\begin{equation}\label{sec:App.++.2}
\E\,\exp\left\{i\,b (y_\zs{t}-y_\zs{s})\right\}
=
\E\,\exp\left\{i\,b (\xi_\zs{t}-\xi_\zs{s})\right\}\frac{\Upsilon_\zs{t}(\xi)}{\Upsilon_\zs{s}(\xi)}
\,.
\end{equation}
To check this equality note that
the process
$$
\Upsilon_\zs{t}(\xi)=\exp
\left\{\int^{t}_\zs{0}\,\frac{S(u)}{\varrho_\zs{1}}\,\d w_\zs{u}
-\,\int^{t}_\zs{0}\,
\frac{S^{2}(u)}{2\varrho^{2}_\zs{1}}\,
\d u
\right\}
$$
is the gaussian martingale. From here we directly obtain the squation \eqref{sec:App.++.2}. Hence
Proposition \ref{Pr.sec:App.1++}.
\endproof

\bigskip

\bigskip

\medskip

\medskip


\end{document}